\documentclass[a4paper,11pt,final]{article} 

\usepackage{graphicx}
\usepackage{amssymb,amsmath,amsthm,mathrsfs,xfrac}
\usepackage{enumitem}
\usepackage{a4wide}
\usepackage{hyperref}
\usepackage{cleveref}
\usepackage{array}
\usepackage{xcolor}

\numberwithin{equation}{section}
\numberwithin{figure}{section}
\allowdisplaybreaks[4]

\newtheorem{proposition}{Proposition}[section]
\newtheorem{theorem}[proposition]{Theorem}
\newtheorem{lemma}[proposition]{Lemma}
\newtheorem{definition}[proposition]{Definition}
\newtheorem{corollary}[proposition]{Corollary}
\newtheorem{remark}[proposition]{Remark}

\newenvironment{proofof}[1]{\smallskip\noindent{\textbf{Proof~of~#1.}}%
  \hspace{1pt}}{\hspace{-5pt}{\nobreak\quad\nobreak\hfill\nobreak%
    $\square$\vspace{2pt}\par}\smallskip\goodbreak}

\newcommand{\C}[1]{\mathbf{C}^{#1}}
\newcommand{\Cc}[1]{\mathbf{C}_c^{#1}}

\newcommand{\BV}{\mathbf{BV}}

\renewcommand{\L}[1]{\mathbf{L}^#1}
\newcommand{\Lloc}[1]{{\mathbf{L}_{\mathbf{loc}}^{#1}}}

\newcommand{\W}[2]{{\mathbf{W}^{#1,#2}}}

\newcommand{\modulo}[1]{{\left|#1\right|}}
\newcommand{\norma}[1]{{\left\|#1\right\|}}
\newcommand{\caratt}[1]{{\chi_{\strut#1}}}
\newcommand{\reali}{{\mathbb{R}}}
\newcommand{\naturali}{{\mathbb{N}}}

\renewcommand{\epsilon}{\varepsilon}
\renewcommand{\phi}{\varphi}
\renewcommand{\theta}{\vartheta}
\renewcommand{\O}{\mathcal{O}(1)}
\newcommand{\tv}{\mathinner{\rm TV}}
\newcommand{\spt}{\mathop{\rm spt}}

\newcommand{\wto}{\rightharpoonup}

\renewcommand{\d}[1]{\mathinner{\mathrm{d}{#1}}}
\renewcommand{\div}{\mathinner{\nabla\cdot}} 
\newcommand{\grad}{\mathinner{\nabla}}

\newcommand{\1}{\mathbf{1}}

\newcommand{\function}[5]{
    \begin{array}{@{}l|ccl@{}}#1: & #2 & \longrightarrow & #3 \\[5pt] & \displaystyle{#4} & \longmapsto & \displaystyle{#5} \end{array}}

\newcommand{\sOmega}{\mathop{\smash[b]{\,*_{{\strut\Omega}\!}}}}
\newcommand{\sOmegaF}{\mathop{\,*_{{\strut\Omega}\!}}}

\DeclareMathOperator*{\esssup}{ess\,sup}

\makeatletter
\let\@fnsymbol\@arabic
\makeatother

\title{Non Local Mixed Systems with Neumann Boundary Conditions}

\author{Rinaldo M.~Colombo\footnotemark[1] \and Elena Rossi\footnotemark[2] \and Abraham Sylla\footnotemark[1]} \date{ }

\begin{document}
\maketitle
\footnotetext[1]{INdAM Unit \& Department of Information Engineering,
  University of Brescia, via Branze, 38, 25123 Brescia, Italy. Email:
  \texttt{rinaldo.colombo@unibs.it} and
  \texttt{abraham.sylla@unibs.it}} \footnotetext[2]{University of
  Modena and Reggio Emilia, INdAM Unit \& Department of Sciences and
  Methods for Engineering, Via Amendola 2 -- Pad.~Morselli, 42122
  Reggio Emilia, Italy. Email: \texttt{elena.rossi13@unimore.it}}

\begin{abstract}
  \noindent We prove well posedness and stability in $\L1$ for a class
  of mixed hyperbolic-parabolic non linear and non local equations in
  a bounded domain with \emph{no flow} along the boundary. While the
  treatment of boundary conditions for the hyperbolic equation is
  standard, the extension to $\L1$ of classical results about
  parabolic equations with Neumann conditions is here achieved.

  \medskip

  \noindent\textit{2020~Mathematics Subject Classification:} 35M30,
  35L04, 35K20

  \medskip


\end{abstract}


\section{Introduction}
\label{sec:Intro}

We consider the following non linear and non local problem on a
bounded domain $\Omega \subset \reali^n$
\begin{equation}
  \label{eq:1}
  \left\{
    \begin{array}{l}
      \partial_t u
      + \div \left(u \,v\left(t,w (t) \right)\right)
      =
      \alpha \left(t,x,w (t)\right) \, u + a(t,x) \,,
      \\
      \partial_t w
      - \mu \, \Delta w
      =
      \beta \left(t,x,u (t),w (t)\right) \, w + b (t,x)\,,
    \end{array}
  \right.
  \quad (t,x) \in [0,T] \times \Omega \,.
\end{equation}
When $n=2,3$, this mixed system is motivated by a variety of
predator-prey models. Indeed, for instance, $u = u (t,x)$ can be the
density of a population that \emph{chases} the other population
$w = w (t,x)$. This chase is described through the non local operator
$v$, able to model the movement of $u$ towards regions where the
concentration of $w$ is higher. A change in the sign of $v$ allows to
model the case where the population $u$ \emph{escapes} from $w$. The
$w$ population diffuses isotropically in all directions. The source
terms $\alpha$ and $\beta$ account for natality, mortality or
predation, while $a$ and $b$ may describe controls acting on the
system, which consist in the introduction of the two species at desired
times and locations.

We equip problem~\eqref{eq:1} with the following initial data and
conditions at the boundary:
\begin{equation}
  \label{eq:ibd}
  \begin{cases}
    u(0,x) = u_o(x)
    \\
    w(0,x)= w_o(x)
  \end{cases}
  \quad x \in \Omega,
  \qquad\quad
  \begin{cases}
    u(t,\xi) = 0
    \\
    \nabla w (t,\xi) \cdot \nu (\xi) =0
  \end{cases}
  \quad (t,\xi) \in \mathopen]0,T\mathclose[ \times \partial\Omega.
\end{equation}
The present choice of Neumann boundary conditions for the parabolic
equation is motivated by the \emph{no flow} requirement typically
suited to the physical setting considered.  Recall that boundary
conditions for the hyperbolic equation bear an entirely different
meaning, since they can be essentially neglected whenever
characteristics exit the domain, see~\cite{MR542510, MR2322018,
  MR3819847, MR1884231}. Thus, the first boundary condition
in~\eqref{eq:ibd} does not prevent $u$ from exiting $\Omega$.

The analytical treatment of~\eqref{eq:1} relies on that of the
separate hyperbolic and parabolic equations
\begin{equation}
  \label{eq:33}
  \partial_t u + \nabla \cdot\left(u \, c (t,x)\right) = A (t,x) \, u + a (t,x)
  \quad \mbox{ and } \quad
  \partial_t w - \mu \, \Delta w = B (t,x) \, w + b (t,x)
\end{equation}
together with the initial and boundary
conditions~\eqref{eq:ibd}. While general well posedness results
related to the former equation are available in the literature, the
treatment in $\L1$ of the latter equation with Neumann boundary
conditions has received far less attention. Indeed, since the
classical books~\cite{FriedmanBook, MR0241822}, the literature offers
a variety of results in $\L2$, a choice which is hardly justifiable in
the present physical setting. Here, on the contrary, the $\L1$ norm
has a clear meaning but $\L1$ stability estimates were not available,
not even in~\cite{FriedmanBook, MR0241822}, especially in the case of
Neumann boundary conditions. Also the general mixed
hyperbolic-parabolic setting in~\cite{MR1880522} does not comprise the
well posedness of~\eqref{eq:1}. Therefore, we provide a definition of
weak solution to the parabolic equation and, correspondingly, we
develop a well posedness and stability theory in the $\L1$ norm. In
this procedure, clearly, the many properties that follow from
reflexivity can not be used. However, it is remarkable that the weak
completeness of $\L1$ plays a key role in our treatment of the
parabolic equation.

Note that in~\eqref{eq:1} the dependence of $v$, $\alpha$ and $\beta$
on $u (t)$ and $w (t)$ is of a \emph{functional} nature, allowing for
\emph{non local} dependencies. Indeed, $u (t)$ and $w (t)$
in~\eqref{eq:1} denote the \emph{functions} $x \mapsto u (t,x)$ and
$x \mapsto w (t,x)$, both defined on all $\Omega$.

Systems of this form arise, for instance, in predator-prey
models~\cite{parahyp} and can be used in the control of parasites,
see~\cite{SAPM2021, Pfab2018489}. A similar mixed
hyperbolic-parabolic system is considered, in one space dimension,
in~\cite{MR4273477}, where Euler equations substitute the balance law
in~\eqref{eq:1}. A mixed ODE -- parabolic PDE predator-prey model with
Neumann boundary conditions is presented in~\cite{MR3556642}: here the
predators' movement is the superposition of a directed hunting and a
random dispersion. The present model~\eqref{eq:1}--\eqref{eq:ibd} is
applicable also to the setting of pursuit-evasion games, similarly
to~\cite{MR4373985}. Here, however, the movement of the pursuer is not
purely diffusive, but it is directed towards the average gradient of
the evaders' density.

A mixed hyperbolic-parabolic system motivated by population dynamics
is considered in~\cite{MR2483341}, where the theoretical framework is
set in $\L2$ and the theory of $m$-accretive operators is the key
analytic tool. A more applied result is~\cite{MR3485306}, where the
description of an aneurysm leads to a mixed hyperbolic-parabolic
system in $1$ space dimension. Global classical solutions to a
parabolic predator-prey system, under Neumann boundary conditions, are
exhibited in~\cite{MR4385435}, motivated by the dynamics of competing
populations with repulsive chemotaxis. In the $\L2$ framework, local
in time well posedness of a mixed hyperbolic-parabolic system is
obtained in~\cite{MR1357364}, by means of a Cauchy sequence of
approximate solutions. A mixed elliptic-parabolic problem also with
biological motivations is studied in~\cite{MR4589866}.

The next section presents the main analytical
results. System~\eqref{eq:1} is split in the $2$
equations~\eqref{eq:33}. The former one is dealt with by means of
results mainly coming from the literature,
see~\S~\ref{sec:hyperbolic-problem}. On the contrary, in
\S~\ref{sec:parabolic-estimates} we present results on the parabolic
part developed \emph{ad hoc} in $\L1$ for the purposes of the present
work. Finally, all proofs are deferred to
Section~\ref{sec:analytical-proofs}.

\section{Main Results}
\label{sec:main}

Throughout, the following notation is used.
$\reali_+ = [0, +\infty\mathclose[$. If $A \subseteq \reali^n$, the
characteristic function $\caratt{A}$ is defined by
$\caratt{A} (x) = 1$ if and only if $x \in A$ and $\caratt{A} (x) = 0$
if and only if $x \in \reali^n\setminus A$. For $x_o \in \reali^n$ and
$r>0$, $B (x_o,r)$ is the open sphere centered at $x_o$ with radius
$r$.

Fix $T > 0$. We pose the following assumptions on $\Omega$ and on the
functions appearing in problem~\eqref{eq:1}:
\begin{enumerate}[label=$\mathbf{(\Omega)}$]
\item\label{it:omega} $\Omega$ is a non empty, bounded and connected
  open subset of $\reali^n$, with $\C{1,\gamma}$ boundary, for a
  $\gamma \in \mathopen]0,1\mathclose]$.
\end{enumerate}
\begin{enumerate}[label=$\mathbf{(\boldsymbol{v})}$]
\item\label{it:v}
  $v\colon [0,T] \times \L1 (\Omega; \reali) \to (\C2 \cap \W1\infty)
  (\Omega; \reali^n)$ is such that for a constant $K_v>0$ and for a
  map $k_v \in \Lloc\infty ([0,T]\times\reali_+; \reali_+)$ non
  decreasing in each argument, for all $t,t_1,t_2 \in [0,T]$ and
  $w,w_1, w_2 \in \L1 (\Omega;\reali)$,
  \begin{eqnarray*}
    \norma{v (t,w)}_{\L\infty (\Omega;\reali^n)}
    & \leq
    & K_v \, \norma{w}_{\L1 (\Omega;\reali)}
    \\
    \norma{D_x v (t,w)}_{\L\infty (\Omega;\reali^n)}
    & \leq
    & K_v\, \norma{w}_{\L1 (\Omega;\reali)}
    \\
    \norma{v (t_1,w_1) - v (t_2,w_2)}_{\L\infty (\Omega; \reali^n)}
    & \leq
    & K_v \left(
      \modulo{t_1 - t_2} + \norma{w_1 - w_2}_{\L1 (\Omega; \reali)}\right)
    \\
    \norma{D^2_x v (t,w)}_{\L1 (\Omega;\reali^{n\times n})}
    & \leq
    & k_v \left(t,\norma{w}_{\L1 (\Omega;\reali)}\right) \,
      \norma{w}_{\L1 (\Omega;\reali)}
    \\
    \norma{\div \left(v (t_1,w_1) - v (t_2,w_2)\right)}_{\L\infty (\Omega; \reali)}
    & \leq
    & k_v \left(t,\max_{i=1,2}\norma{w_i}_{\L1 (\Omega;\reali)}\right) \,
      \norma{w_1 - w_2}_{\L1 (\Omega;\reali)} \,.
  \end{eqnarray*}
\end{enumerate}
\begin{enumerate}[label=$\mathbf{(\boldsymbol{\alpha})}$]
\item\label{it:alpha}
  $\alpha \colon [0,T] \times \Omega \times \W11(\Omega;\reali) \to
  \reali$ admits a constant $K_\alpha> 0$ such that, for
  a.e.~$t \in [0,T]$, for a.e.~$x \in \Omega$ and all
  $w_1,w_2 \in \W11 (\Omega;\reali)$,
  \begin{displaymath}
    \sup_{x \in \Omega} \modulo{\alpha (t,x,w_1) - \alpha (t,x,w_2)}
    \leq
    K_\alpha \; \norma{w_1-w_2}_{\L1 (\Omega;\reali)}
  \end{displaymath}
  and there exists $k_\alpha \in \L1 ([0,T];\reali_+)$ such that for
  all $w \in \W11 (\Omega;\reali)$ and a.e.~$t \in [0,T]$,
  \begin{eqnarray*}
    \tv\left(\alpha (t, \cdot, w)\right)
    & \leq
    & K_\alpha \left(1 + \norma{\nabla w}_{\L1 (\Omega;\reali^n)}\right)
    \\
    \esssup_{x\in\Omega} \modulo{\alpha (t,x,w)}
    & \leq
    & k_\alpha (t) \left(1 + \norma{w}_{\L1 (\Omega;\reali)}\right) \,.
  \end{eqnarray*}

\end{enumerate}
\begin{enumerate}[label=$\mathbf{(\boldsymbol{a})}$]
\item\label{it:a}
  $a \in \L1 \left([0,T]; \L\infty(\Omega; \reali)\right)$ and
  $t \mapsto \tv\left(a (t)\right)$ is in $\L1 ([0,T];\reali)$.
\end{enumerate}
\begin{enumerate}[label=$\mathbf{(\boldsymbol{\beta})}$]
\item\label{it:betaFinal}
  $\beta \colon [0,T] \times \Omega \times \L1 (\Omega;\reali)\times
  \L1 (\Omega;\reali) \to \reali$ admits a constant $K_\beta > 0$ and
  a $k_\beta \in \L1 ([0,T];\reali_+)$ such that, for
  a.e.~$t \in [0,T]$, all $x \in \Omega$ and
  $u, u_1,u_2,w, w_1,w_2 \in \L1 (\Omega;\reali)$,
  \begin{eqnarray*}
    \modulo{\beta (t,x,u_1,w_1) - \beta (t,x,u_2,w_2)}
    & \leq
    & K_\beta \left(
      \norma{u_1-u_2}_{\L1 (\Omega;\reali)}
      +
      \norma{w_1-w_2}_{\L1 (\Omega;\reali)}
      \right)
    \\
    \esssup_{x\in\Omega} \modulo{\beta (t,x,u,w)}
    & \leq
    & k_\beta (t) \,.
  \end{eqnarray*}

\end{enumerate}
\begin{enumerate}[label=$\mathbf{(\boldsymbol{b})}$]
\item\label{it:b} $b \in \L1 ([0,T] \times \Omega;\reali)$.
\end{enumerate}

\noindent Since~\eqref{eq:1} is the coupling of a hyperbolic and a
parabolic problem, the following definition of solution
to~\eqref{eq:1} is a sort of gluing of the definitions of solutions to
the equations in~\eqref{eq:33}

\begin{definition}
  \label{def:def_sol}
  A pair $(u,w) \in \C0\left([0,T]; \L1 (\Omega;\reali^2)\right)$ is a
  solution to problem~\eqref{eq:1}--\eqref{eq:ibd} if, setting
  \begin{displaymath}
    c(t,x) = v\left(t, w(t)\right) (x),
    \qquad
    A(t,x) = \alpha\left(t,x,w(t)\right),
    \qquad
    B(t,x) =
    \beta\left(t,x,u(t),w(t)\right) \,,
  \end{displaymath}
  the function $u$ solves, according to Definition~\ref{def:hypGen},
  the problem
  \begin{equation}
    \label{eq:21}
    \left\{
      \begin{array}{l@{\qquad}l}
        \partial_t u + \nabla \cdot \left(u \; c (t,x)\right) = A (t,x) u + a(t,x)
        &  (t,x) \in [0,T] \times \Omega
        \\
        u (t,\xi) = 0
        & (t,\xi) \in \mathopen]0,T\mathclose[ \times\partial\Omega
        \\
        u (0,x) = u_o (x)
        & x \in \Omega
      \end{array}
    \right.
  \end{equation}
  and the function $w$ solves, according to
  Definition~\ref{def:paraGen}, the problem
  \begin{equation}
    \label{eq:13}
    \left\{
      \begin{array}{l@{\qquad}l}
        \partial_t w - \mu \; \Delta w = B (t,x) w + b(t,x)
        & (t,x) \in [0,T] \times \Omega
        \\
        \nabla w (t,\xi) \cdot \nu (\xi) = 0
        & (t,\xi) \in \mathopen]0,T\mathclose[ \times\partial\Omega
        \\
        w (0,x) = w_o (x)
        & x \in \Omega.
      \end{array}
    \right.
  \end{equation}
\end{definition}

We now state the main result of this work: it ensures the well
posedness in $\L1$ of~\eqref{eq:1} and provides stability estimates to
be used, for instance, in control problems based on~\eqref{eq:1}.

\begin{theorem}
  \label{thm:main}
  Fix $T > 0$.  Assume that~\ref{it:omega}, \ref{it:v},
  \ref{it:alpha}, \ref{it:a}, \ref{it:betaFinal} and~\ref{it:b}
  hold. Let
  $(u_o,w_o) \in (\L1\cap \BV) (\Omega;\reali) \times \L1
  (\Omega;\reali)$. Then:
  \begin{enumerate}[label=\bf(M\arabic*)]
  \item \label{item:27} Problem~\eqref{eq:1} admits a unique solution
    $(u,w)$ in the sense of Definition~\ref{def:def_sol}.

  \item \label{item:31} The map $t \mapsto (u,w) (t)$ is in
    $\C0 \left([0,T]; \L1 (\Omega;\reali^2)\right)$.

  \item \label{item:28} If
    $(u_o^1,w_o^1), (u_o^2,w_o^2) \in (\L1\cap \BV) (\Omega;\reali)
    \times \L1 (\Omega;\reali)$, the corresponding solutions
    $(u_1,w_1)$ and $(u_2,w_2)$ satisfy the estimate:
    \begin{displaymath}
      \norma{(u_1,w_1) (t) - (u_2,w_2) (t)}_{\L1 (\Omega;\reali^2)}
      \leq
      C (t) \; \norma{(u_o^1,w_o^1) - (u_o^2,w_o^2)}_{\L1 (\Omega;\reali^2)}
    \end{displaymath}
    where $C \in \L\infty ([0,T];\reali_+)$ depends on $\mu$,
    \ref{it:omega}, \ref{it:v}, \ref{it:alpha}, \ref{it:betaFinal},
    \ref{it:a}, \ref{it:b} and on $\tv (u_o^1)$, $\tv (u_o^2)$,
    $\norma{(u_o^1,w_o^1)}_{\L1 (\Omega;\reali^2)}$,
    $\norma{(u_o^2,w_o^2)}_{\L1 (\Omega;\reali^2)}$.
  \item \label{item:29} If $a_1,a_2$ satisfy~\ref{it:a} and $b_1,b_2$
    satisfy~\ref{it:b}, the corresponding solutions $(u_1,w_1)$ and
    $(u_2,w_2)$ satisfy the estimate:
    \begin{displaymath}
      \norma{(u_1,w_1) (t) {-} (u_2,w_2) (t)}_{\L1 (\Omega;\reali^2)}
      \leq
      C (t) \,
      \left(
        \norma{a_1{-}a_2}_{\L1 ([0,t]\times\Omega; \reali)}
        +
        \norma{b_1{-}b_2}_{\L1 ([0,t]\times\Omega; \reali)}
      \right)
    \end{displaymath}
    where $C \in \L\infty ([0,T];\reali_+)$ depends on $\mu$,
    \ref{it:omega}, \ref{it:v}, \ref{it:alpha}, \ref{it:betaFinal},
    \ref{it:a}, \ref{it:b} and on $\tv (u_o)$,
    $\norma{(u_o,w_o)}_{\L1 (\Omega;\reali^2)}$.

  \item \label{item:33} If $k_\beta$ in~\ref{it:betaFinal} is bounded
    and $u_o \geq 0$, $w_o \geq 0$, $a\geq 0$ and $b \geq 0$, then for
    a.e.~$t \in [0,T]$, the solution $(u,w)$ satisfies $u (t) \geq 0$
    and $w (t) \geq 0$.
  \end{enumerate}
\end{theorem}

\noindent Remark that if the various assumptions hold on the time
interval $\reali_+$, then Theorem~\ref{thm:main} ensures well
posedness on $\reali_+$, with the function $C$ appearing
in~\ref{item:28} and~\ref{item:29} belonging to
$\Lloc\infty (\reali_+, \reali_+)$.

To deal with non local operators on a bounded domain $\Omega$, the
modified convolution introduced in~\cite[\S~3]{siam2018} is an
adequate tool. For a function $\rho \in \L1 (\Omega;\reali)$ and a
smooth kernel $\eta$, it reads
\begin{equation}
  \label{eq:58}
  (\rho \sOmega \eta) (x)
  =
  \dfrac{\int_\Omega \rho (y) \; \eta (x-y) \d{y}}{\int_{\Omega} \eta (x-y) \d{y}} \,.
\end{equation}
The quantity $(\rho \sOmegaF \eta) (x)$ is an average of the values
attained by $\rho$ in $\Omega$ around $x$ as soon as the kernel $\eta$
satisfies, for instance,
\begin{description}
\item[($\boldsymbol{\eta}$)] $\eta (x) = \tilde\eta (\norma{x})$,
  where $\tilde\eta \in \C3 (\reali_+; \reali)$,
  $\spt \tilde\eta = [0, \ell]$, $\ell > 0$, $\tilde\eta' \leq 0$,
  $\tilde\eta' (0) = \tilde\eta'' (0) = 0$ and
  $\int_{\reali^n} \eta (\xi) \d\xi = 1$.
\end{description}

\noindent It is often reasonable to assume that $u$ hunts $w$ moving
towards areas with higher density of $w$, or else that $u$ escapes
from $w$ towards regions with lower $w$ density. Thus, $v$ is parallel
to the average gradient of $w$ in $\Omega$, such as
\begin{equation}
  \label{eq:61}
  v (t,w)
  = k (t) \;
  \dfrac{\grad (w \sOmegaF \eta)}{\sqrt{1 + \norma{\grad (w \sOmega \eta)}^2}}
\end{equation}
where, for instance,
\begin{equation}
  \label{eq:31}
  \eta (x) = \overline{\eta} \; \left(1 -
    (\norma{x}/ \ell)^4\right)^4 \; \caratt{B (0, \ell)} (x) \,.
\end{equation}
Here, $\ell$ has the clear physical meaning of the distance, or
\emph{horizon}, at which individuals of the $u$ population \emph{feel}
the presence of the $w$ population.  The normalization parameter
$\overline\eta$ is chosen so that
$\int_{\reali^n} \eta (x) \d{x} = 1$.  A choice like~\eqref{eq:61} is
consistent with the requirements~\ref{it:v}, as proved
in~\cite[Lemma~3.2]{siam2018}.

\subsection{The Hyperbolic Problem}
\label{sec:hyperbolic-problem}

We focus on the hyperbolic problem~\eqref{eq:21}. For completeness, we
present the standard definition of solution and a detailed well
posedness result based on the current literature. Precise references
are provided in~\S~\ref{subsec:proofs-relat-hyperb}.

\begin{definition}
  \label{def:hypGen}
  A function $u \in \L\infty ([0,T] \times \Omega; \reali)$ is a
  solution to~\eqref{eq:21} if for any test function
  $\phi \in \Cc1 (\mathopen] -\infty, T\mathclose[ \times \Omega;
  \reali)$,
  \begin{equation}
    \label{eq:23}
    \int_0^T \int_\Omega
    \left(u \,
      (\partial_t \phi  + c  \cdot \nabla \phi )
      +
      (A \, u + a )\phi
    \right)
    \d{x}\d{t}
    + \int_\Omega u_o (x) \, \phi(0,x) \d{x}
    = 0.
  \end{equation}
\end{definition}

\begin{theorem}
  \label{thm:hypGen}
  Fix $T>0$. Assume~\ref{it:a} and
  \begin{enumerate}[label={\bf(c)}]
  \item\label{it:H1}
    $c \in \left(\C0 \cap \L\infty\right) ([0,T] \times
    \Omega;\reali^n)$, $c(t) \in \C1 (\Omega; \reali^n)$ for all
    $t \in [0,T]$,
    $D_x c \in \L\infty ([0,T] \times \Omega; \reali^{n\times n})$.
  \end{enumerate}
  \begin{enumerate}[label={\bf(A)}]
  \item\label{it:H2}
    $A \in \L1 \left([0,T]; \L\infty(\Omega;\reali)\right)$ and
    $t \mapsto \tv\left(A (t)\right)$ is in $\L1 ([0,T];\reali)$.
  \end{enumerate}
  Then, for all $u_o \in \L1 (\Omega;\reali)$, problem~\eqref{eq:21}
  admits a unique solution in the sense of
  Definition~\ref{def:hypGen}. Moreover
  \begin{enumerate}[label=\bf(H\arabic*)]
  \item \label{item:11} For all $t \in [0,T]$,
    \begin{equation}
      \label{eq:17}
      \norma{u (t)}_{\L1 (\Omega;\reali)}
      \leq
      \left(
        \norma{u_o}_{\L1 (\Omega;\reali)}
        +
        \norma{a}_{\L1 ([0,t]\times\Omega; \reali)}
      \right) \,
      \exp\left(\norma{A}_{\L1 ([0,t];\L\infty(\Omega;\reali))}\right) \,.
    \end{equation}
  \item \label{item:20} For all $t \in [0,T]$, if
    $u_o \in \L\infty (\Omega;\reali)$,
    \begin{equation}
      \label{eq:25}
      \begin{array}{rcl}
        \norma{u (t)}_{\L\infty (\Omega;\reali)}
        & \leq
        & \left(
          \norma{u_o}_{\L\infty (\Omega;\reali)}
          +
          \norma{a}_{\L1 ([0,t]; \L\infty(\Omega; \reali))}
          \right)
        \\
        &
        & \quad \times
          \exp\left(
          \norma{A}_{\L1 ([0,t];\L\infty(\Omega;\reali))}
          + \norma{\nabla\cdot c}_{\L1 ([0,t];\L\infty (\Omega;\reali))}
          \right) \,.
      \end{array}
    \end{equation}
  \item \label{item:12} If $u_o \geq 0$ and $a \geq 0$, then
    $u \geq 0$.
  \item \label{item:13} If $u_o^1,u_o^2 \in \L1 (\Omega;\reali)$, then
    the corresponding solutions $u_1,u_2$ satisfy for all
    $t \in [0,T]$,
    \begin{equation}
      \norma{u_1 (t) - u_2 (t)}_{\L1 (\Omega;\reali)}
      \leq
      \norma{u_o^1 - u_o^2}_{\L1 (\Omega;\reali)}
      \exp\left(\norma{A}_{\L1 ([0,t]; \L\infty(\Omega;\reali))}\right)\,.
    \end{equation}
  \item \label{item:14} If $A_1,A_2$ satisfy~\ref{it:H2} and $a_1,a_2$
    satisfy~\ref{it:a}, the corresponding solutions $u_1$ and $u_2$
    satisfy for all $t \in [0,T]$
    \begin{eqnarray*}
      &
      & \norma{u_1 (t) - u_2 (t)}_{\L1 (\Omega;\reali)}
      \\
      & \leq
      & \exp \left(
        \max
        \left\{
        \norma{A_1}_{\L1 ([0,t]; \L\infty(\Omega;\reali))}
        \,,\;
        \norma{A_2}_{\L1 ([0,t]; \L\infty(\Omega;\reali))}
        \right\}\right)
      \\
      &
      & \quad \times
        \left(
        \norma{u_o}_{\L1 (\Omega;\reali)}
        +
        \norma{a_1}_{\L1 ([0,t]; \L\infty(\Omega;\reali))}
        \right)
        \norma{A_2-A_1}_{\L1 ([0,t];\L\infty(\Omega;\reali))}
      \\
      &
      & + \exp\left(\norma{A_2}_{\L1 ([0,t];\L\infty(\Omega;\reali))}\right)
        \norma{a_1 - a_2}_{\L1 ([0,t]\times\Omega;\reali)} \,.
    \end{eqnarray*}

  \item \label{item:16} If $u_o \in \BV (\Omega;\reali)$,
    $c (t) \in \C2 (\Omega;\reali^n)$ for all $t \in [0,T]$ and
    $\nabla \nabla \cdot c \in \L1 ([0,T]\times\Omega;\reali^n)$, then
    the map $t \to u (t)$ is locally Lipschitz continuous in
    $\L1 (\Omega;\reali)$.

  \item \label{item:15} If $u_o \in \BV (\Omega;\reali)$, $c_1,c_2$
    satisfy~\ref{it:H1} and
    $c_1 (t),\, c_2 (t) \in \C2 (\Omega;\reali^n)$ for all
    $t \in [0,T]$ and
    $\nabla \nabla \cdot c_1, \nabla \nabla \cdot c_2 \in \L1
    ([0,T]\times\Omega;\reali^n)$, then the corresponding solutions
    $u_1, \, u_2$ satisfy
    \begin{align*}
      \!\!\!
      & \norma{u_2 (t) - u_1 (t)}_{\L1 (\Omega;\reali)}
      \\
      \!\!\!
      \leq
      & C \!\left(\!
        \norma{D_x c_1}_{\L1 ([0,t];\L\infty (\Omega; \reali^{n\times n}))},
        \norma{\nabla \nabla\cdot c_1}_{\L1 ([0,t]\times\Omega;\reali^n)},
        \norma{A}_{\L1 ([0,t];\L\infty (\Omega;\reali))},
        \!\int_0^t \!\!\!\tv \! A (\tau)\d\tau \!
        \! \right)
      \\
      \!\!\!
      & \times\left(
        \norma{c_1 - c_2}_{\L1 ([0,t]; \L\infty (\Omega;\reali^n))}
        +
        \norma{\div(c_1 - c_2)}_{\L1 ([0,t]; \L\infty (\Omega;\reali))}
        \right)
    \end{align*}
    where $C$ also depends on
    $\norma{a}_{\L1 ([0,t];\L\infty(\Omega;\reali))}$,
    $\int_0^t \tv\!\left(a (\tau)\right)\d\tau$ and
    $\norma{u_o}_{\L1 (\Omega;\reali)}$,
    $\norma{u_o}_{\L\infty (\Omega;\reali)}$, $\tv (u_o)$.
  \end{enumerate}
\end{theorem}

\noindent The proof is based on results from the literature detailed
in~\S~\ref{subsec:proofs-relat-hyperb}.

\subsection{The Parabolic Problem}
\label{sec:parabolic-estimates}

We focus on the $\L1$ well posedness for the parabolic
problem~\eqref{eq:13}, first adapting the classical definition of weak
solution, see for instance~\cite[\S~2.3]{MR3017032}, to the case of
interest here.

\begin{definition}
  \label{def:paraGen}
  A function $w \in \L\infty \left([0,T] ; \L1(\Omega;\reali) \right)$
  with $w (t) \in \W11 (\Omega;\reali)$ for a.e.~$t \in [0,T]$ is a
  \emph{weak solution to~\eqref{eq:13}} if for all test functions
  $\phi \in \Cc\infty ([0,T\mathclose[ \times \overline\Omega;\reali)$
  \begin{equation}
    \label{eq:22}
    \begin{array}{rcl}
      \displaystyle
      \int_0^T \int_\Omega w \; \partial_t \phi \d{x} \d{t}
      - \mu \int_0^T \int_\Omega \nabla w \cdot \nabla \phi \d{x} \d{t}
      \\
      \displaystyle
      +\int_0^T \int_\Omega B \, w \; \phi \d{x} \d{t}
      +\int_0^T \int_\Omega b \; \phi \d{x} \d{t}
      +\int_\Omega w_o (x)\; \phi (0,x) \d{x}
      & =
      & 0 \,.
    \end{array}
  \end{equation}
\end{definition}

A relevant consequence of the definition chosen above is the following
convergence result where the weak completeness of $\L1$ is
essential. This Lemma is of use in a few key points in the proof of
Theorem~\ref{thm:paraGen}.

\begin{lemma}
  \label{lem:weakComplete}
  For $h \in \naturali$, let $b_h \in \L1 ([0,T]\times\Omega;\reali)$,
  $B_h \in \L1\left([0,T];\L\infty (\Omega;\reali)\right)$ and
  $w_o^h \in \L1 (\Omega;\reali)$ be such that~\eqref{eq:13} admits
  the solution $w_h$ in the sense of
  Definition~\ref{def:paraGen}. Moreover, assume that
  \begin{displaymath}
    \begin{array}{r@{\,}c@{\,}ll@{\qquad}r@{\,}c@{\,}ll}
      \lim_{h\to+\infty} b_h
      & =
      & b
      & \mbox{ in } \L1 ([0,T]\times\Omega;\reali) \,,
      & \lim_{h\to+\infty} B_h
      & =
      & B
      & \mbox{ in } \L1\left([0,T];\L\infty (\Omega;\reali)\right) \,,
      \\
      \lim_{h\to+\infty} w_o^h
      & =
      & w_o
      & \mbox{ in }\L1 (\Omega;\reali) \,,
      & \lim_{h\to+\infty} w_h
      & =
      & w
      & \mbox{ in } \L\infty\left([0,T];\L1 (\Omega;\reali)\right)\,.
    \end{array}
  \end{displaymath}
  Then:
  \begin{enumerate}[label=\bf(\arabic*)]
  \item \label{item:18} For a.e.~$t\in [0,T]$,
    $w (t) \in \W11 (\Omega;\reali)$.
  \item \label{item:25} $\lim_{h\to+\infty} \nabla w_h = \nabla w$
    weakly in $\L1 ([0,T]\times\Omega;\reali)$.
  \item \label{item:26} $w$ is a solution to~\eqref{eq:13} in the
    sense of Definition~\ref{def:paraGen}.
  \end{enumerate}
\end{lemma}

\noindent The following result differs from others found in the
literature in its being set in $\L1$ and in its referring to Neumann
boundary conditions.

\begin{theorem}
  \label{thm:paraGen}
  Let $\Omega$ satisfy~\ref{it:omega}. Fix $T,\mu>0$. Assume
  $B \in \L1 \left([0,T]; \L\infty(\Omega;\reali)\right)$ and
  $b \in \L1([0,T]\times\Omega;\reali)$. Then, for all
  $w_o \in \L1 (\Omega;\reali)$, problem~\eqref{eq:13} admits a unique
  solution in the sense of Definition~\ref{def:paraGen}. Moreover:
  \begin{enumerate}[label=\bf(P\arabic*)]

  \item \label{item:17} It also holds that
    $w \in \C0 \left([0,T];\L1 (\Omega;\reali)\right)$,
    $\nabla w \in \L1 ([0,T]\times\Omega;\reali^n)$ and
    \begin{equation}
      \label{eq:29}
      \begin{array}{rcl}
        \norma{\nabla w}_{\L1 ([0,T]\times\Omega;\reali^n)}
        & \leq
        & \displaystyle
          \frac{1}{\mu}
          \norma{B}_{\L1 ([0,T];\L\infty (\Omega;\reali))} \;
          \norma{w}_{\L\infty ([0,T];\L1 (\Omega;\reali))}
        \\
        &
        & \displaystyle
          + \frac{1}{\mu} \, \norma{b}_{\L1 ([0,T]\times\Omega;\reali)}
          + \frac{1}{\mu} \norma{w_o}_{\L1 (\Omega;\reali)}  .
      \end{array}
    \end{equation}

  \item \label{item:7} The following implicit representation formula
    holds:
    \begin{equation}
      \label{eq:14}
      \begin{array}{rcl}
        w (t,x)
        & =
        & \displaystyle
          \int_{\Omega} N (t,x,0,y) \; w_o (y) \d{y}
        \\
        &
        & \displaystyle +
          \int_0^t \int_\Omega N (t,x,s,y) \;
          \left(B (s,y) \; w (s,y) + b (s,y)\right) \d{y} \d{s}
      \end{array}
    \end{equation}
    and the Neumann function $N$ is defined in
    Proposition~\ref{prop:Green}.
  \item \label{item:APBound} There exists a positive $K$ depending on
    $\mu$, $\Omega$ --- hence on $n$ --- such that the following
    \emph{a priori} bound holds for all $t \in [0, T]$:
    \begin{equation}
      \label{eq:APBound}
      \!\!\!\!\!\!\!
      \norma{w (t)}_{\L1 (\Omega;\reali)}
      \leq
      K \left(\norma{w_o}_{\L1 (\Omega;\reali)} +
        \norma{b}_{\L1 ([0, t] \times \Omega;\reali)} \right)
      \exp\left(K \, \norma{B}_{\L1 ([0, t]; \L\infty (\Omega;\reali))}\right) \,.
    \end{equation}

  \item \label{item:9} If $w_o^1, w_o^2 \in \L1 (\Omega;\reali)$, the
    corresponding solutions $w_1$ and $w_2$ satisfy for
    all~$t \in [0,T]$,
    \begin{equation}
      \label{eq:15}
      \norma{w_1 (t) - w_2 (t)}_{\L1 (\Omega;\reali)}
      \leq
      K \, \norma{w_o^1-w_o^2}_{\L1 (\Omega;\reali)}\,
      \exp\left(K \, \norma{B}_{\L1([0, t]; \L\infty (\Omega;\reali))} \right) \,.
    \end{equation}
  \item \label{item:10} If
    $B_1, B_2 \in \L1 \left([0,T]; \L\infty(\Omega;\reali)\right)$,
    the corresponding solutions $w_1$ and $w_2$ satisfy for
    all~$t \in [0,T]$,
    \begin{eqnarray}
      \nonumber
      &
      & \norma{w_1 (t) - w_2 (t)}_{\L1 (\Omega;\reali)}
      \\
      \label{eq:16}
      & \leq
      & K^2
        \exp\left(
        K
        (\norma{B_1}_{\L1 ([0, t]; \L\infty (\Omega;\reali))}
        +
        \norma{B_2}_{\L1 ([0, t]; \L\infty (\Omega;\reali))})
        \right)
      \\
      \nonumber
      &
      & \times
        \left(
        \norma{w_o}_{\L1 (\Omega;\reali)}
        +
        \norma{b}_{\L1 ([0,t]\times\Omega;\reali)}
        \right)
        \norma{B_1- B_2}_{\L1([0,t];\L\infty  (\Omega;\reali))} \,.
    \end{eqnarray}

  \item \label{item:32} If
    $b_1, b_2 \in \L1([0,T]\times\Omega;\reali)$, the corresponding
    solutions $w_1$, $w_2$ satisfy for all~$t \in [0,T]$,
    \begin{equation}
      \label{eq:30}
      \norma{w_1 (t) - w_2 (t)}_{\L1 (\Omega;\reali)} \leq
      K
      \exp\left(K \, \norma{B}_{\L1 ([0, t]; \L\infty
          (\Omega;\reali))} \right)
      \norma{b_1 - b_2}_{\L1 ([0,t]\times\Omega;\reali)} \,.
    \end{equation}

  \item \label{item:8} Assume
    $B \in \L\infty ([0,T]\times\Omega;\reali)$. If $b \geq 0$ and
    $w_o \geq 0$, then $w \geq 0$.
  \end{enumerate}
\end{theorem}

\noindent The proof is deferred to \S~\ref{subsec:proofs-relat-parab}.

\section{Analytical Proofs}
\label{sec:analytical-proofs}

\subsection{Hyperbolic Problem}
\label{subsec:proofs-relat-hyperb}

\begin{proofof}{Theorem~\ref{thm:hypGen}}
  The existence and uniqueness of $u$ follow
  from~\cite[Proposition~3.9]{ElenaDirichlet}.

  The \emph{a priori} $\L1$ and $\L\infty$ bounds~\ref{item:11}
  and~\ref{item:20} are obtained in~\cite[Lemma~4.2]{teo}.  Positivity
  in~\ref{item:12} is proved as in~\cite[Lemma~3.12]{ElenaDirichlet}.
  The Lipschitz continuous dependence on the initial
  datum~\ref{item:13} follows from~\ref{item:11} by linearity.  The
  stability estimate~\ref{item:14} is proved through the same
  computations as in~\cite[Lemma~4.3]{teo}, taking advantage of
  linearity and of~\ref{item:11}. The continuity~\ref{item:16} follows
  from~\cite[Lemma~3.13]{ElenaDirichlet}. The stability
  in~\ref{item:15} is proved in~\cite[Lemma~3.1]{ElenaDirichlet}.
\end{proofof}

\subsection{Parabolic Problem}
\label{subsec:proofs-relat-parab}

\begin{proofof}{Lemma~\ref{lem:weakComplete}}
  Using Definition~\ref{def:paraGen}, pass to the limit
  $h \to + \infty$ in
  \begin{eqnarray*}
    \mu \int_0^T \int_\Omega \nabla w_h \cdot \nabla \phi \d{x} \d{t}
    & =
    & \int_0^T \int_\Omega w_h \; \partial_t \phi \d{x} \d{t}
      + \int_0^T \int_\Omega B_h \; w_h \; \phi \d{x} \d{t}
    \\
    &
    & + \int_0^T \int_\Omega b_h \; \phi \d{x} \d{t} +\int_\Omega
      w_o^h(x) \; \phi (0,x) \d{x} \,.
  \end{eqnarray*}
  By the Dominated Convergence Theorem, for all test functions
  $\phi \in \Cc\infty ([0,T\mathclose[ \times \overline\Omega;\reali)$
  \begin{flalign*}
    \int_0^T \int_\Omega w_h \; \partial_t \phi \d{x} \d{t} %
    & \underset{h\to+\infty}{=} \int_0^T \int_\Omega w \; \partial_t
    \phi \d{x} \d{t}%
    & \mbox{[Since $w_h \to w$ in $\L1$]}
    \\
    \int_0^T \int_\Omega B_h \; w_h \; \phi \d{x} \d{t}%
    & \underset{h\to+\infty}{=} \int_0^T \int_\Omega B \; w \; \phi
    \d{x} \d{t} %
    & \mbox{[Since $w_h \to w$ and $B_h \to B$]}
    \\
    \int_0^T \int_\Omega b_h \; \phi \d{x} \d{t} %
    & \underset{h\to+\infty}{=} \int_0^T \int_\Omega b \; \phi \d{x}
    \d{t} %
    & \mbox{[Since $b_h \to b$ in $\L1$]}
    \\
    \int_\Omega w_o^h(x) \; \phi (0,x) \d{x} %
    &\underset{h\to+\infty}{=} \int_\Omega w_o(x) \; \phi (0,x)
    \d{x} %
    & \mbox{[Since $w_o^h \to w_o$ in $\L1$]}
  \end{flalign*}
  As a consequence, we obtain that
  \begin{displaymath}
    \lim_{h\to+\infty}
    \int_0^T \int_\Omega \nabla w_h \cdot \nabla \phi \d{x} \d{t}
    =
    c_\phi
  \end{displaymath}
  for a real number $c_\phi$, so that by the weak completeness of
  $\L{1}$, see~\cite[Corollary~14]{MR1144277}, there exists a map
  $z \in \L1 ([0,T]\times\Omega; \reali^n)$ such that
  $\nabla w_h \underset{h\to+\infty}{\wto} z$ in $\L1$. Choose now
  $\psi \in \Cc\infty(\mathopen]0,T\mathclose[ \times \Omega;\reali)$:
  \begin{flalign*}
    \int_0^T \int_\Omega w \; \nabla \psi \d{x} \d{t} & =
    \lim_{h\to+\infty} \int_0^T \int_\Omega w_h \; \nabla \psi \d{x}
    \d{t} & \mbox{[Since $w_h \underset{h\to+\infty}{\wto} w$ in
      $\L1$]}
    \\
    & = - \lim_{h\to+\infty} \int_0^T \int_\Omega \nabla w_h \; \psi
    \d{x} \d{t} & \mbox{[Since $w_h \in \W11$]}
    \\
    & = - \int_0^T \int_\Omega z \; \psi \d{x} \d{t} & \mbox{[Since
      $\nabla w_h \underset{h\to+\infty}{\wto} z$ in $\L1$]}
  \end{flalign*}
  proving that $w (t) \in \W11 (\Omega;\reali)$ for a.e.~$t \in [0,T]$
  and $\nabla w \in \L1 ([0,T]\times\Omega; \reali^n)$.  Hence, $w$
  satisfies the regularity requirements in
  Definition~\ref{def:paraGen}. Since it also satisfies~\eqref{eq:22},
  $w$ is a weak solution to~\eqref{eq:13} in the sense of
  Definition~\ref{def:paraGen}.
\end{proofof}

We first consider problem~\eqref{eq:13} with $B=0$, namely
\begin{equation}
  \label{eq:2}
  \left\{
    \begin{array}{l@{\qquad}l}
      \partial_t w - \mu \; \Delta w = b (t,x)
      & (t,x) \in [0,T] \times \Omega
      \\
      \nabla w (t,\xi) \cdot \nu (\xi) = 0
      & (t,\xi) \in \mathopen]0,T\mathclose[ \times\partial\Omega
      \\
      w (0,x) = w_o (x)
      & x \in \Omega
    \end{array}
  \right.
\end{equation}
under condition~\ref{it:omega} on $\Omega$.

\begin{lemma}
  \label{lem:weakCont}
  Fix $w_o \in \L1 (\Omega;\reali)$,
  $b \in \L1 ([0,T]\times\Omega;\reali)$ and let $w$
  solve~\eqref{eq:2} in the sense of
  Definition~\ref{def:paraGen}. Then, for all
  $\eta \in \C\infty(\overline\Omega;\reali)$
  \begin{equation}
    \label{eq:10}
    \begin{array}{r@{\,}c@{\,}l@{\qquad}l}
      \forall \, t
      & \in
      & \mathopen]0,T\mathclose[:
      & \displaystyle
        \lim_{h\to 0} \int_\Omega
        \left(w (t+h,x) - w (t,x)\right) \; \eta (x) \d{x}
        = 0
      \\
      t
      & =
      & 0:
      & \displaystyle
        \lim_{h\to 0+} \int_\Omega
        \left(w (h,x) - w_o (x)\right) \; \eta (x) \d{x}
        = 0\,.
    \end{array}
  \end{equation}
\end{lemma}

\begin{proofof}{Lemma~\ref{lem:weakCont}}
  Fix $t \in [0,T\mathclose[$ and introduce the sequence of functions
  \begin{displaymath}
    \chi_k \in \Cc\infty (\reali; [0,1]) \,,\quad
    \spt \chi_k \subseteq [-1,t]
    \quad \mbox{ and } \quad
    \chi_k \underset{k\to+\infty}{\to} \1_{[0, t]}
    \mbox{ pointwise a.e. on } [0,T]\,.
  \end{displaymath}
  Then use~\eqref{eq:22} with $B\equiv 0$, first
  $\phi (s,x) = \chi_k (s+h) \, \eta (x)$ and then
  $\phi (s,x) = \chi_k (s) \, \eta (x)$, for a suitable
  $\eta \in \C\infty(\overline\Omega;\reali)$, with
  $t \in [0,T\mathclose[$ and $h$ sufficiently small. Taking the
  difference of the resulting expressions, we have:
  \begin{eqnarray*}
    &
    &\int_0^T \int_\Omega w (s,x)
      \left(\partial_t \chi_k (s+h) - \partial_t \chi_k (s)\right)
      \eta (x) \d{x} \d{s}
    \\
    & =
    & \mu \int_0^T \int_\Omega \left(\chi_k (s+h) - \chi_k (s)\right)
      \nabla w (s,x) \cdot \nabla \eta (x) \d{x} \d{s}
    \\
    &
    & -
      \int_0^T \int_\Omega b (s,x)
      \left(\chi_k (s+h) - \chi_k (s)\right) \eta (x) \d{x} \d{s}
      +
      \int_\Omega w_o (x) \; \left(\chi_k (0) - \chi_k (h)\right) \; \eta (x)
      \d{x}\,.
  \end{eqnarray*}
  If $t>0$, the latter term above vanishes and in the limit
  $k \to +\infty$ the first equality in~\eqref{eq:10} follows.

  When $t=0$, the above terms reduce to
  \begin{eqnarray*}
    &
    &\int_0^T \int_\Omega w (s,x) \;
      \partial_t \chi_k (s+h) \;
      \eta (x) \d{x} \d{s}
    \\
    & =
    & \mu \int_0^T \int_\Omega \chi_k (s+h) \;
      \nabla w (s,x) \cdot \nabla \eta (x) \d{x} \d{s}
    \\
    &
    & -
      \int_0^T \int_\Omega b (s,x)
      \; \chi_k (s+h) \; \eta (x) \d{x} \d{s}
      -
      \int_\Omega w_o (x) \; \chi_k (h) \; \eta (x) \d{x}
  \end{eqnarray*}
  and as $k \to +\infty$ the second equality in~\eqref{eq:10} follows.
\end{proofof}

When the initial datum $w_o$ is in $\L2 (\Omega;\reali)$ and the
source $b$ is in $\L2 ([0,T]\times\Omega;\reali)$, strong $\L2$
continuity in time is available.

\begin{lemma}
  \label{lem:salsa}
  If $b \in \L2 \left([0,T];\L2 (\Omega;\reali)\right)$ and
  $w_o \in \L2 (\Omega;\reali)$, then problem~\eqref{eq:2} admits a
  unique solution $w \in \L2 \left([0,T];\W12 (\Omega;\reali)\right)$
  in the sense of Definition~\ref{def:paraGen} and
  $\partial_t w \in \L2 \left([0,T];\W12 (\Omega;\reali)^*\right)$,
  so that $w \in \C0\left([0,T];\L2 (\Omega;\reali)\right)$.
\end{lemma}

\begin{proofof}{Lemma~\ref{lem:salsa}}
  Use~\cite[Problem~(10.35)]{MR3362185} and apply~\cite[point~b) in
  Theorem~7.104]{MR3362185}. Note that the definition of weak solution
  in~\cite[p.~592]{MR3362185} implies Definition~\ref{def:paraGen} due
  to the density of $\Cc\infty (\Omega;\reali)$ in
  $\W12 (\Omega;\reali)$ and by the Dominated Convergence Theorem.
\end{proofof}

We introduce for later use the space $\mathcal{V}$ as the closure of
$\W12 (\reali\times\Omega;\reali)$ with respect to the norm
\begin{displaymath}
  \norma{w}_\mathcal{V} =
  \norma{\nabla w}_{\L2 (\reali\times\Omega;\reali)}
  +
  \esssup_{t \in \reali} \norma{w (t)}_{\L2 (\Omega;\reali)} \,.
\end{displaymath}

\begin{proposition}
  \label{prop:Green}
  Let $\Omega$ satisfy~\ref{it:omega} and $\mu>0$. Then, there exists
  a function
  \begin{displaymath}
    N \in \C0\left( \left\{(t,x,s,y) \in (\reali \times \Omega)^2
        \colon (t,x) \neq (s,y)\right\}; \reali \right)
  \end{displaymath}
  such that
  \begin{enumerate}[label=\bf(N\arabic*)]

  \item \label{item:2} For all $(s,y) \in \reali\times\Omega$,
    \begin{displaymath}
      (t,x) \mapsto N (t,x,s,y)
      \in
      \Lloc1 (\reali \times \Omega; \reali)
      \mbox{ and }
      (t,x) \mapsto \nabla_x N (t,x,s,y)
      \in
      \Lloc1 (\reali\times\Omega;\reali^n)\,.
    \end{displaymath}

  \item \label{item:6} There exist positive $C$, $\kappa$ and $K$ such
    that for all $(t,x,s,y) \in (\reali_+ \times\Omega)^2$
    with $t>s$
    \begin{eqnarray}
      \label{eq:32}
      \modulo{N (t,x,s,y)}
      & \leq
      & C \; \left(1+\dfrac{1}{(t-s)^{n/2}}\right) \,
        \exp\left(-\dfrac{\kappa\norma{x-y}^2}{t-s}\right) \,;
      \\
      \label{eq:5}
      \int_\Omega \modulo{N (t,x,s,y)} \d{x}
      & \leq
      & K.
    \end{eqnarray}

  \item \label{item:5} Fix $T>0$. For every
    $w_o \in \L2 (\Omega;\reali)$ and all
    $b \in \C\infty ([0,T]\times\overline\Omega;\reali)$, the map
    \begin{equation}
      \label{eq:3}
      w (t,x)
      =
      \int_{\Omega} N (t,x,0,y) \; w_o (y) \d{y}
      +
      \int_0^t \int_\Omega N (t,x,s,y) \; b (s,y) \d{y} \d{s}
    \end{equation}
    is the unique function in $\mathcal{V}$ satisfying~\eqref{eq:22},
    with $B \equiv 0$, for all test function
    $\phi \in \Cc\infty ([0,T\mathclose[ \times
    \overline\Omega;\reali)$.
  \end{enumerate}
\end{proposition}

\begin{proofof}{Proposition~\ref{prop:Green}}
  Here we apply~\cite[Theorem~3.9]{MR3017032}, thus we need to verify
  conditions~(A1)--(A2) in~\cite{MR3017032}. (A1) holds by~\cite[(1)
  in Examples~4.1.1]{MR3017032} thanks
  to~\cite[Theorem~9.7]{BrezisBook} which can be applied thanks
  to~\ref{it:omega}; (A2) holds by~\cite[(1) in
  Examples~4.1.2]{MR3017032}, that applies in the present scalar
  case. We can apply~\cite[Theorem~3.21]{MR3017032} since also~(A3)
  in~\cite{MR3017032} holds by~\cite[(3) in
  Examples~4.1.3]{MR3017032}, since the coefficients in~\eqref{eq:2}
  are constant and~\ref{it:omega} holds.

  The regularity of $N$ and the proof of~\ref{item:2} directly follow
  from~\cite[Theorem~3.9]{MR3017032}. To prove~\eqref{eq:32}, start
  from the last line in the proof of~\cite[Theorem~3.21]{MR3017032}:
  \begin{displaymath}
    \modulo{N (t,x,s,y)}
    \leq
    C \;\max \left\{1, \frac{t-s}{{\rm diam} (\Omega)}\right\}^{n/2} \;
    (t-s)^{-n/2} \;
    \exp\left(-\dfrac{\kappa\norma{x-y}^2}{t-s}\right)
  \end{displaymath}
  which implies~\eqref{eq:32}, up to relabeling $C$. Now, \eqref{eq:5}
  is a direct consequence with
  $K = C \, \left({\rm meas} (\Omega) +(\pi/\kappa)^{n/2}\right)$.

  Consider~\ref{item:5}. The proof that the expression~\eqref{eq:3}
  solves~\eqref{eq:2} in the sense of Definition~\ref{def:paraGen}
  follows from~\cite[Theorem~3.9]{MR3017032}, as well as uniqueness.

\end{proofof}

\begin{remark}
  \label{rem:pos}
  {\rm A close look at the proof of~\cite[Theorem~3.9]{MR3017032}
    shows also that $N\geq0$. This positivity is not explicitly
    considered in~\cite{MR3017032} since the object of that work is a
    system of parabolic equations and thus $N$ results to be matrix
    valued.  This property, though basic, is not necessary in the
    sequel, hence we omit its proof whose rigorous exposition might
    significantly lengthen this work. However, $N\geq0$ ensures that
    \begin{displaymath}
      w_o \in \L2 (\Omega;\reali_+)\quad \mbox{ and } \quad
      b \in \C\infty ([0,T]\times\overline\Omega;\reali_+)
      \implies w\geq 0 \,.
    \end{displaymath}
  }
\end{remark}

\begin{proposition}
  \label{prop:L2toL1}
  Let $\Omega$ satisfy~\ref{it:omega}. Fix $T,\mu>0$. For every
  $w_o \in \L1 (\Omega;\reali)$ and any
  $b \in \L1([0, T] \times \Omega ; \reali)$, problem~\eqref{eq:2}
  admits a unique weak solution $w$ in the sense of
  Definition~\ref{def:paraGen}. Moreover,
  \begin{enumerate}[label=\bf(\arabic*)]

  \item \label{item:22} $w$ admits the representation~\eqref{eq:3}
    with $N$ as defined in Proposition~\ref{prop:Green}.

  \item \label{item:23}
    $w \in \C0\left([0,T];\L1 (\Omega;\reali)\right)$ and
    $\nabla w \in \L1 ([0,T]\times\Omega;\reali^n)$.
  \end{enumerate}
\end{proposition}

\begin{proofof}{Proposition~\ref{prop:L2toL1}}
  We split the proof in a few steps.

  \paragraph{Step~1: Uniqueness.} Assume $w_1,w_2$ solve~\eqref{eq:2}
  in the sense of Definition~\ref{def:paraGen}. Then, their difference
  $w$ satisfies, for all test function
  $\phi \in \Cc\infty ([0,T[\times \overline\Omega;\reali)$,
  \begin{equation}
    \label{eq:WeakUniqueness}
    \int_0^T \int_\Omega w \; \partial_t \phi \d{x} \d{t}
    -\mu \int_0^T \int_\Omega \nabla w \cdot \nabla \phi \d{x} \d{t}=0.
  \end{equation}
  First apply~\eqref{eq:WeakUniqueness} with a test function depending
  only on time to obtain, thanks to Lemma~\ref{lem:weakCont}, that for
  all $t \in [0, T]$,
  \begin{equation}
    \label{eq:8}
    \int_{\Omega} w(t, x) \d{x} = 0 \,.
  \end{equation}
  Then, fix arbitrary $a,b \in \mathopen]0,T\mathclose[$ with
  $a<b$. Introduce a sequence of test functions
  \begin{displaymath}
    \chi_h \in \Cc\infty (\mathopen]0,T\mathclose[; [0,1]) \,,\quad
    \spt \chi_h \subseteq [a,b]
    \quad \mbox{ and } \quad
    \chi_h \to \1_{[a, b]} \mbox{ pointwise a.e. on } [0,T]\,.
  \end{displaymath}
  Let $\psi \in \C{\infty}(\overline\Omega; \reali)$ and
  apply~\eqref{eq:WeakUniqueness} with $\phi = \chi_h \, \psi$,
  obtaining
  \begin{eqnarray*}
    0
    & =
    & \int_a^b \int_\Omega w \; \partial_t \chi_h \; \psi \d{x} \d{t}
      -
      \mu \int_a^b \int_\Omega \chi_h \; \nabla w \cdot \nabla \psi \d{x} \d{t}
    \\
    & \underset{h\to+\infty}{=}
    & \int_\Omega w (b,x) \d{x} - \int_\Omega w (a,x) \d{x}
      - \mu \int_a^b \int_\Omega  \nabla w \cdot \nabla \psi \d{x} \d{t}
    \\
    & =
    & - \mu \int_a^b \int_\Omega  \nabla w \cdot \nabla \psi \d{x} \d{t}
  \end{eqnarray*}
  where we used~\eqref{eq:8}. Hence, there exists a $c \in \reali^n$
  such that for a.e.~$x\in\Omega$,
  \begin{equation}
    \label{eq:9}
    \int_a^b \nabla w (t,x) \d{t} = c
  \end{equation}
  and, by the arbitrariness of $a$ and $b$, it must be $c=0$.
  From~\eqref{eq:9} we thus obtain that for a.e.~$x \in \Omega$ and
  for all $g \in \L\infty ([0,T];\reali^n)$
  \begin{displaymath}
    \int_0^T \nabla w (t,x) \cdot g (t) \d{t} = 0
  \end{displaymath}
  so that $\nabla w (t,x) =0$ for all $t \in [0,T]$ and
  a.e.~$x \in \Omega$. Together with~\eqref{eq:8} and the
  connectedness of $\Omega$, this implies that $w=0$, proving
  uniqueness.

  \paragraph{Step~2: Approximation.}
  By~\cite[(1.8) in \S~1.14]{Giusti}, there exists a sequence
  $b_h \in \Cc\infty (\reali^{1+n}; \reali)$ such that $b_h \to b$ in
  $\L1 ([0,T]\times\Omega;\reali)$ as $h\to+\infty$. Similarly, there
  exists a sequence $w_o^h \in \Cc\infty (\reali^n, \reali)$ such that
  $w_o^h \to w_o$ in $\L1 (\Omega;\reali)$ as $h\to+\infty$. Note that
  $w_o^h \in \L2 (\Omega;\reali)$ for all $h$.

  For all $h \in \naturali$, by~\ref{item:5} in
  Proposition~\ref{prop:Green} there is a unique solution
  $w_h \in \mathcal{V}$ to \eqref{eq:2} with source $b_h$ and initial
  datum $w_o^h$, which is given by the representation
  \begin{equation}
    \label{eq:11}
    w_h(t,x) = \int_{\Omega} N (t,x,0,y) \; w_o^h (y) \d{y} + \int_0^t
    \int_\Omega N (t,x,s,y) \; b_h(s,y) \d{y} \d{s} \,.
  \end{equation}
  Then, for $h,k \in \naturali$ with $k >h$, thanks to~\eqref{eq:5} we
  have
  \begin{eqnarray*}
    \norma{w_k (t) - w_h (t)}_{\L1 (\Omega;\reali)}
    & \leq
    & \int_\Omega \int_{\Omega} \modulo{N (t,x,0,y)} \; \modulo{w_o^k (y) - w_o^h (y)} \d{y} \d{x}
    \\
    &
    & + \int_\Omega  \int_0^t
      \int_\Omega \modulo{N (t,x,s,y)} \; \modulo{b_k(s,y) - b_h(s,y)}
      \d{y} \d{s} \d{x}
    \\
    & =
    & \int_\Omega \left(\int_{\Omega} \modulo{N (t,x,0,y)} \d{x}\right) \;
      \modulo{w_o^k (y) - w_o^h (y)} \d{y}
    \\
    &
    & + \int_\Omega \int_0^t
      \left(\int_\Omega \modulo{N (t,x,s,y)} \d{x}\right)
      \modulo{b_k(s,y) - b_h(s,y)}  \d{s} \d{y}
    \\
    & \leq
    & K
      \left(
      \norma{w_o^k - w_o^h}_{\L1 (\Omega;\reali)}
      +
      \norma{b_k - b_h}_{\L1 ([0,T]\times\Omega;\reali)}
      \right)
  \end{eqnarray*}
  proving that there exists a function
  $w \in \L\infty\left([0,T];\L1 (\Omega;\reali)\right)$ such that
  \begin{equation}
    \label{eq:24}
    \lim_{h\to+\infty} \norma{w_h - w}_{\L\infty\left([0,T];\L1 (\Omega;\reali)\right)} =0 \,.
  \end{equation}

  \paragraph{Step~3: Existence in
    $\L\infty \left([0,T];\L1 (\Omega;\reali)\right)$.} To prove that
  $w$ solves~\eqref{eq:2} in the sense of Definition~\ref{def:paraGen}
  it is now sufficient to apply Lemma~\ref{lem:weakComplete} with
  $B_h\equiv 0$.

  Hence, $w$ satisfies the regularity requirements in
  Definition~\ref{def:paraGen}. Since it also satisfies~\eqref{eq:22},
  $w$ is a weak solution to~\eqref{eq:2} in the sense of
  Definition~\ref{def:paraGen}, with $B \equiv 0$.

  \paragraph{Step~4: $\L1(\Omega;\reali)$ Continuity in Time.} Note
  that the sequence $w_h$ defined above also satisfies
  $w_h \in \L2 \left([0,T];\W12 (\Omega;\reali)\right)$ and
  $\partial_t w_h \in \L2 \left([0,T];\W12
    (\Omega;\reali)^*\right)$, so that
  $w_h \in \C0\left([0,T];\L2 (\Omega;\reali)\right)$ by
  Lemma~\ref{lem:salsa}. Recall that
  $\C0\left([0,T];\L2 (\Omega;\reali)\right) \subseteq
  \C0\left([0,T];\L1 (\Omega;\reali)\right)$, hence the
  convergence~\eqref{eq:24} then ensures that
  $w \in \C0\left([0,T];\L1 (\Omega;\reali)\right)$.

  \paragraph{Step~5: Representation Formula.}
  It is now sufficient to pass to the limit $h\to+\infty$
  in~\eqref{eq:11} to prove that $w$ admits the
  representation~\eqref{eq:3}.

  The proof is completed.
\end{proofof}

\begin{corollary}
  \label{cor:properties}
  Let~\ref{it:omega} hold. Fix $T,\mu>0$. Let
  $w_o^1, w_o^2 \in \L1 (\Omega;\reali)$ and
  $b_1, b_2 \in \L1([0, T] \times \Omega ; \reali)$. Call $w_1, w_2$
  the corresponding solutions to~\eqref{eq:2}. Then, for
  all~$t \in [0,T]$,
  \begin{equation}
    \label{eq:12}
    \norma{w_1 (t) - w_2 (t)}_{\L1 (\Omega;\reali)}
    \leq
    K
    \left(
      \norma{w_o^1 - w_o^2}_{\L1 (\Omega;\reali)}
      +
      \norma{b_1 - b_2}_{\L1 ([0,t]\times\Omega;\reali)}
    \right)
  \end{equation}
  where $K$ is as in~\ref{item:6} of Proposition~\ref{prop:Green}.
\end{corollary}

\noindent Simply apply~\ref{item:22} of Proposition~\ref{prop:L2toL1}
to the difference $w_2-w_1$ and exploit the linearity of~\eqref{eq:2}.

\medskip

\begin{proofof}{Theorem~\ref{thm:paraGen}}
  We split the proof in a few steps.

  \paragraph*{Step 1: Problem~\eqref{eq:13} admits a unique solution on
    $[0,T]$ satisfying~\ref{item:17} and~\ref{item:7}.}
  Consider the operators
  \begin{equation}
    \label{eq:20}
    \function{\Lambda}{\C0\left([0,T];\L1(\Omega;\reali)\right)}{
      \L1([0,T] \times \Omega;\reali)}{w}{\Lambda \, w \quad
      \mbox{ where }(\Lambda \, w)(t, x) = B(t, x) \, w(t, x) + b (t,x)}
  \end{equation}
  and, with reference to Definition~\ref{def:paraGen},
  \begin{displaymath}
    \function{\Phi}{\L1([0,T] \times \Omega;\reali)}{
      \C0\left([0,T];\L1(\Omega;\reali)\right)}{\beta}{w \quad
      \mbox{ where }   \left\{
        \begin{array}{@{}l@{\qquad}l@{}}
          \partial_t w - \mu \; \Delta w = \beta (t,x)
          & (t,x) \in [0,T] \times \Omega
          \\
          \nabla w (t,\xi) \cdot \nu (\xi) = 0
          & (t,\xi) \in \mathopen]0,T\mathclose[ \times\partial\Omega
          \\
          w (0,x) = w_o (x)
          & x \in \Omega \,.
        \end{array}
      \right.}
  \end{displaymath}
  Let us precise that $\Lambda$ is well defined, meaning that if
  $w \in \C0\left([0,T];\L1(\Omega;\reali)\right)$ then, by the
  assumptions on $B$ and $b$,
  $\Lambda(w) \in \L1([0,T] \times \Omega;\reali)$. Similarly, $\Phi$
  is well defined by Proposition~\ref{prop:L2toL1}.

  By the assumption on $B$, there exist times $t_i$ for
  $i=0, \ldots, m$ such that
  $0=t_0<t_1<\cdots <t_i < t_{i+1} < \cdots < t_m =T$ and
  \begin{equation}
    \label{eq:19}
    \int_{t_i}^{t_{i+1}} \norma{B (s)}_{\L\infty (\Omega;\reali)} \d{s}
    <
    \frac1{2 \, K}
  \end{equation}
  with $K$ as in~\eqref{eq:5} in Proposition~\ref{prop:Green}.
  Clearly, thanks to the continuity proved in
  Proposition~\ref{prop:L2toL1}, $w$ in
  $\L\infty \left([0,T];\L1 (\Omega;\reali)\right)$ with
  $w (t) \in \W11 (\Omega;\reali)$ for a.e.~$t \in [0,T]$ is a weak solution
  to~\eqref{eq:13} if and only if it is a fixed point of
  $\Phi \circ \Lambda$ in
  $\C0 \left([0,T];\L1 (\Omega;\reali)\right)$. To construct such
  fixed point we apply Banach Fixed Point Theorem iteratively in each
  of the spaces $\C0 \left([t_i,t_{i+1}]; \L1 (\Omega;\reali)\right)$
  for $i = 0, \ldots, m-1$.

  Notice that for all
  $w_1, w_2 \in \C0\left([t_i,t_{i+1}];\L1(\Omega;\reali)\right)$ with
  $w_1 (t_i) = w_2 (t_i)$, for all $t \in [t_i,t_{i+1}]$,
  \begin{flalign*}
    & \norma{\Phi \circ \Lambda(w_1)(t) - \Phi \circ
      \Lambda(w_2)(t)}_{\L1(\Omega;\reali)}
    \\
    &\leq K \int_{t_i}^t \norma{\Lambda(w_1)(s) -
      \Lambda(w_2)(s)}_{\L1(\Omega;\reali)} \; \d{s} %
    &[\mbox{By Corollary~\ref{cor:properties}}]
    \\
    & = K \int_{t_i}^t \int_{\Omega} \modulo{B(s, y)} \;
    \modulo{w_1(s, y) - w_2(s, y)} \, \d{y} \d{s} %
    &[\mbox{By~\eqref{eq:20}}]
    \\
    &\leq K \int_{t_i}^t \norma{B (s)}_{\L\infty (\Omega;\reali)}
    \int_{\Omega} \modulo{w_1(s, y) - w_2(s, y)} \; \d{y} \d{s}
    \\
    &\leq K \int_{t_i}^t \norma{B (s)}_{\L\infty (\Omega;\reali)} \;
    \d{s} \; \norma{w_1 - w_2}_{\C0
      ([t_i,t_{i+1}];\L1(\Omega;\reali))}
    \\
    &\leq \frac12 \; \norma{w_1 - w_2}_{\C0
      ([t_i,t_{i+1}];\L1(\Omega;\reali))} \,.%
    & [\mbox{By~\eqref{eq:19}]}
  \end{flalign*}
  An iterated application of Banach Fixed Point Theorem ensures the
  existence of $w_* \in \C0\left([0, T] ; \L1(\Omega;\reali)\right)$
  such that $w_* = \Phi \circ \Lambda(w_*)$.

  Define $\tilde b = B \, w_* + b$, so that
  $\tilde b \in \L1 ([0,T]\times\Omega;\reali)$. By construction,
  $w_*$ solves
  \begin{displaymath}
    \left\{
      \begin{array}{l@{\qquad}l}
        \partial_t w - \mu \; \Delta w = \tilde b (t,x)
        & (t,x) \in [0,T] \times \Omega
        \\
        \nabla w (t,\xi) \cdot \nu (\xi) = 0
        & (t,\xi) \in \mathopen]0,T\mathclose[ \times\partial\Omega
        \\
        w (0,x) = w_o (x)
        & x \in \Omega
      \end{array}
    \right.
  \end{displaymath}
  in the sense of Definition~\ref{def:paraGen}. Hence,
  $w_* (t) \in \W11 (\Omega;\reali)$ for a.e.~$t \in [0,T]$ and
  $\nabla w \in \L1 ([0,T]\times\Omega;\reali^n)$ by~\ref{item:23} in
  Proposition~\ref{prop:L2toL1}.  This shows that $w_*$ is a weak
  solution to \eqref{eq:13} in the sense of
  Definition~\ref{def:paraGen} on $[0,T]$ and~\ref{item:17} holds.

  To prove the bound~\eqref{eq:29}, by~\eqref{eq:22} we have
  \begin{eqnarray*}
    \!\!\!
    &
    & \modulo{
      \int_0^T \int_\Omega
      [w \quad -\mu\, \nabla w] \Bigl[
      \begin{array}{@{}c@{}}
        {\partial_t \phi}\\{\nabla \phi}
      \end{array}
    \Bigr]
    \d{x} \d{t}
    }
    \\
    \!\!\!
    & \leq
    & \modulo{\int_0^T \int_\Omega (B\, w +b) \phi \d{x}\d{t}}
      + \modulo{\int_\Omega w_o \, \phi (0,\cdot) \d{x}}
    \\
    \!\!\!
    & \leq
    & \left(\!
      \norma{B}_{\L1 ([0,T];\L\infty (\Omega;\reali))} \,
      \norma{w}_{\L\infty ([0,T];\L1 (\Omega;\reali))}
      {+} \norma{b}_{\L1 ([0,T]\times\Omega;\reali)}
      {+} \norma{w_o}_{\L1 (\Omega;\reali)}
      \! \right) \!
      \norma{\phi}_{\L\infty ([0,T]\times\Omega;\reali)}
  \end{eqnarray*}
  which together with
  \begin{displaymath}
    \norma{\nabla w}_{\L1 ([0,T]\times\Omega;\reali^n)}
    \leq
    \frac{1}{\mu} \, \norma{[w \quad -\mu\, \nabla w]}_{\L1 ([0,T]\times\Omega;\reali^{n+1})}
  \end{displaymath}
  completes the proof of~\eqref{eq:29}.

  By construction, \ref{item:7} holds. Step~1 is proved.

  \paragraph*{Step~2: \ref{item:APBound} holds.} Using \eqref{eq:14},
  for all $t \in [0, T]$, by~\eqref{eq:5}
  \begin{displaymath}
    \begin{aligned}
      \norma{w(t)}_{\L1 (\Omega;\reali)} %
      & \leq K \left( \norma{w_o}_{\L1 (\Omega;\reali)} + \int_{0}^t
        \int_{\Omega} \left(\modulo{B(s, y)} \; \modulo{w(s, y)} +
          \modulo{b(s, y)}\right) \; \d{y}\d{s}\right)
      \\
      & \leq K \left( \norma{w_o}_{\L1 (\Omega;\reali)} +
        \norma{b}_{\L1([0, t]; \L1(\Omega;\reali))} +\int_{0}^t
        \norma{B(s)}_{\L\infty (\Omega;\reali)} \; \norma{w(s)}_{\L1
          (\Omega;\reali)} \d{s} \right).
    \end{aligned}
  \end{displaymath}
  An application of Gronwall Lemma leads to~\eqref{eq:APBound}.

  \paragraph*{Step~3: \ref{item:9} holds.} This is an immediate
  consequence of~\ref{item:APBound} by linearity.

  \paragraph{Step~4: \ref{item:10} and~\ref{item:32} hold.}

  Using~\eqref{eq:14}, compute:
  \begin{eqnarray*}
    &
    & \norma{w_1 (t) - w_2 (t)}_{\L1 (\Omega;\reali)}
    \\
    & \leq
    & K \!
      \int_0^t \! \int_\Omega \!
      \modulo{B_1 (s,y) \, w_1 (s,y) - B_2 (s,y)\, w_2 (s,y)}
      \d{y} \d{s}
      {+} K \!
      \int_0^t \! \int_\Omega \!
      \modulo{b_1 (s,y) - b_2 (s,y)}
      \d{y} \d{s}
    \\
    & \leq
    & K
      \int_0^t \int_\Omega \modulo{B_1 (s,y) \, w_1 (s,y) - B_1 (s,y)\, w_2 (s,y)}
      \d{y} \d{s}
    \\
    &
    & + K
      \int_0^t \! \int_\Omega \! \modulo{B_1 (s,y) \, w_2 (s,y) - B_2 (s,y)\, w_2 (s,y)}
      \d{y} \d{s}
      {+} K
      \int_0^t \! \int_\Omega \! \modulo{b_1 (s,y) - b_2 (s,y)}
      \d{y} \d{s}
    \\
    & \leq
    & K \,
      \int_0^t \norma{B_1 (s)}_{\L\infty (\Omega;\reali)} \,
      \norma{w_1 (s) - w_2 (s)}_{\L1 (\Omega;\reali)} \, \d{s}
    \\
    &
    & + K \,
      \norma{B_1- B_2}_{\L1([0,t];\L\infty (\Omega;\reali))} \,
      \norma{w_2}_{\C0([0,t];\L1 (\Omega;\reali))}
      + K \,
      \norma{b_1 - b_2}_{\L1 ([0,t]\times\Omega;\reali)} \,.
  \end{eqnarray*}
  By Gronwall Lemma and using~\ref{item:APBound},
  \begin{eqnarray*}
    &
    & \norma{w_1 (t) - w_2 (t)}_{\L1 (\Omega;\reali)}
    \\
    & \leq
    & K
      \left(
      \norma{B_1- B_2}_{\L1([0,t];\L\infty (\Omega;\reali))} \,
      \norma{w_2}_{\C0([0,t];\L1 (\Omega;\reali))}
      +
      \norma{b_1 - b_2}_{\L1 ([0,t]\times\Omega;\reali)}
      \right)
    \\
    &
    &\times \exp\left(K \,
      \norma{B_1}_{\L1 ([0, t]; \L\infty (\Omega;\reali))}\right)
    \\
    & \leq
    & K \, \norma{b_1 - b_2}_{\L1 ([0,t]\times\Omega;\reali)}
      \exp\left(K \, \norma{B_1}_{\L1 ([0, t]; \L\infty (\Omega;\reali))}\right)
    \\
    &
    & + K^2 \, \norma{B_1- B_2}_{\L1([0,t];\L\infty (\Omega;\reali))}
      \left(\norma{w_o}_{\L1 (\Omega;\reali)}
      + \norma{b_2}_{\L1 ([0,t]\times\Omega;\reali)} \right) \\
    &
    & \times \exp\left(K \left(\norma{B_1}_{\L1 ([0, t]; \L\infty (\Omega;\reali))}
      + \norma{B_2}_{\L1 ([0, t]; \L\infty (\Omega;\reali))}\right) \right),
  \end{eqnarray*}
  which implies~\eqref{eq:16} and~\eqref{eq:30}.

  \paragraph*{Step~5: \ref{item:8} holds with
    $B \in \L\infty ([0,T]\times\Omega;\reali)$.}  By~\cite[(1.8) in
  \S~1.14]{Giusti}, there exists a sequence
  $b_h \in \Cc\infty (\reali^{1+n}; \reali_+)$ such that $b_h \to b$
  in $\L1 ([0,T]\times\Omega;\reali)$ as $h\to+\infty$. Similarly,
  there exists a sequence $w_o^h \in \Cc\infty (\reali^n, \reali_+)$
  such that $w_o^h \to w_o$ in $\L1 (\Omega;\reali)$ as
  $h\to+\infty$. Note that $w_o^h \in \L2 (\Omega;\reali)$ for all
  $h$. Call $w_h$ the solution to
  \begin{displaymath}
    \left\{
      \begin{array}{l@{\qquad}l}
        \partial_t w_h - \mu \; \Delta w_h = B (t,x) w_h + b_h(t,x)
        & (t,x) \in [0,T] \times \Omega
        \\
        \nabla w_h (t,\xi) \cdot \nu (\xi) = 0
        & (t,\xi) \in \mathopen]0,T\mathclose[ \times\partial\Omega
        \\
        w (0,x) = w_o^h (x)
        & x \in \Omega.
      \end{array}
    \right.
  \end{displaymath}
  in the sense of the definition given
  in~\cite[p.~592]{MR3362185}. By~\cite[Remark~10.19]{MR3362185},
  $w_h \geq 0$.

  Using~\ref{item:9} and~\ref{item:32}, we have
  \begin{eqnarray*}
    \norma{w_h (t) -w_k (t)}_{\L1(\Omega;\reali)}
    & \leq
    & K
      e^{K \norma{B}_{\L1([0,t];\L\infty (\Omega;\reali))}}
      \left(
      \norma{w_o^h {-} w_o^k}_{\L1 (\Omega;\reali)}
      {+}
      \norma{b_h {-} b_k}_{\L1 ([0,t]\times\Omega;\reali)}
      \right)
  \end{eqnarray*}
  hence there exists a
  $w \in \L\infty\left([0,T];\L1 (\Omega;\reali)\right)$ such that
  $w_h \to w$ in $\L\infty\left([0,T];\L1 (\Omega;\reali)\right)$.

  Apply now Lemma~\ref{lem:weakComplete} to prove that $w$ is a
  solution to~\eqref{eq:13} in the sense of
  Definition~\ref{def:paraGen}. By construction, $w \geq 0$.

  The proof is completed.
\end{proofof}

\subsection{Mixed Problem}
\label{subsec:proofs-related-mixed}

\begin{proofof}{Theorem~\ref{thm:main}}
  Define $u_1$ and $w_1$ as solutions to
  \begin{displaymath}
    \left\{
      \begin{array}{@{}l@{\quad}r@{\,}c@{\,}l}
        \partial_t u_{1}
        =
        a(t,x)
        & (t,x)
        & \in
        & [0,T] \times \Omega
        \\
        u_1 (t,\xi) = 0
        & (t,\xi)
        & \in
        & \mathopen]0,T\mathclose[ \times \partial \Omega
        \\
        u_1 (0,x) = u_o (x)
        & x
        & \in
        & \Omega
      \end{array}
    \right.
    \qquad
    \left\{
      \begin{array}{l@{\quad}r@{\,}c@{\,}l@{}}
        \partial_t w_{1}
        - \mu \, \Delta w_{1}
        =
        b (t,x)
        & (t,x)
        & \in
        & [0,T] \times \Omega
        \\
        \nabla w_1 (t,\xi) \cdot \nu (\xi) = 0
        & (t,\xi)
        & \in
        & \mathopen]0,T\mathclose[ \times \partial \Omega
        \\
        w_1 (0,x) = w_o (x)
        & x
        & \in
        & \Omega
      \end{array}
    \right.
  \end{displaymath}
  in the sense of Definition~\ref{def:hypGen} and of
  Definition~\ref{def:paraGen}. Clearly, $u_1$ and $w_1$ are well
  defined by Theorem~\ref{def:hypGen} and
  Theorem~\ref{thm:paraGen}. Define recursively, for
  $i \in \naturali\setminus\{0\}$, $u_{i+1}$ as solution to
  \begin{equation}
    \label{eq:6}
    \left\{
      \begin{array}{l@{\qquad}r@{\,}c@{\,}l}
        \partial_t u_{i+1}
        + \div \left(u_{i+1} \, c_i(t,x) \right)
        =
        A_i (t,x)\, u_{i+1} + a(t,x)
        & (t,x)
        & \in
        & [0,T] \times \Omega
        \\
        u_{i+1} (t,\xi) = 0
        & (t,\xi)
        & \in
        & \mathopen]0,T\mathclose[ \times \partial \Omega
        \\
        u_{i+1} (0,x) = u_o (x)
        & x
        & \in
        & \Omega
      \end{array}
    \right.
  \end{equation}
  in the sense of Definition~\ref{def:hypGen} and $w_{i+1}$ as
  solution to
  \begin{equation}
    \label{eq:7}
    \left\{
      \begin{array}{l@{\qquad}r@{\,}c@{\,}l}
        \partial_t w_{i+1}
        - \mu \, \Delta w_{i+1}
        =
        B_i (t,x) \, w_{i+1} + b (t,x)
        & (t,x)
        & \in
        & [0,T] \times \Omega
        \\
        \nabla w_{i+1} (t,\xi) \cdot \nu (\xi) = 0
        & (t,\xi)
        & \in
        & \mathopen]0,T\mathclose[ \times \partial \Omega
        \\
        w_{i+1} (0,x) = w_o (x)
        & x
        & \in
        & \Omega
      \end{array}
    \right.
  \end{equation}
  in the sense of Definition~\ref{def:paraGen}, where for
  $i \in \naturali\setminus\{0\}$,
  \begin{equation}
    \label{eq:18}
    \begin{array}{rcl}
      c_i (t,x)
      & =
      & v \left(t,w_i (t)\right) (x)
      \\
      A_i (t,x)
      & =
      & \alpha \left(t,x,w_i (t)\right)
    \end{array}
    \qquad\qquad
    B_i (t,x)
    =
    \beta \left(t,x,u_i (t),w_i (t)\right) \,.
  \end{equation}
  We aim to prove that $(u_i,w_i)$ are well defined for all
  $i \in \naturali\setminus\{0,1\}$ and constitute a Cauchy sequence
  with respect to the $\C0 \left([0,T] ; \L1 (\Omega;\reali^2)\right)$
  distance as soon as $T$ is sufficiently small.

  \paragraph{Step~1: $c_i$ satisfies~\ref{it:H1} in
    Theorem~\ref{thm:hypGen}.}

  We first check the continuity. By~\eqref{eq:18} and~\ref{it:v} we
  have
  \begin{eqnarray*}
    &
    & \norma{c_i (t,x) - c_i (t_o,x_o)}
    \\
    & \leq
    & \norma{c_i (t,x) - c_i (t_o,x)} + \norma{c_i (t_o,x) - c_i (t_o,x_o)}
    \\
    & =
    & \norma{v\left(t, w_i (t)\right) (x) - v\left(t_o, w_i (t_o)\right) (x)}
      +
      \norma{v\left(t_o, w_i (t_o)\right) (x) - v\left(t_o, w_i (t_o)\right) (x_o)}
    \\
    & \leq
    & K_v \left(\modulo{t-t_o} + \norma{w_i (t) - w_i (t_o)}_{\L1 (\Omega;\reali)}\right)
      +
      \norma{v \left(t_o,w_i (t_o)\right) (x) - v \left(t_o,w_i (t_o)\right) (x_o)}
  \end{eqnarray*}
  where the terms $\norma{w_i (t) - w_i (t_o)}_{\L1 (\Omega;\reali)}$
  and
  $\norma{v \left(t_o,w_i (t_o)\right) (x_o) - v \left(t_o,w_i
      (t_o)\right) (x_o)}$ vanish as $t\to t_o$ and $x \to x_o$ by the
  continuity of $w_i$ and that of $v\left(t_o,w_i (t_o)\right)$,
  ensured by~\ref{it:v}.

  To prove the $\L\infty$ boundedness, observe that by~\ref{it:v}
  \begin{equation}
    \label{eq:28}
    \norma{c_i (t,x)}
    \leq
    \norma{v\left(t, w_i (t)\right)}_{\L\infty (\Omega;\reali^n)}
    \leq
    K_v \; \norma{w_i (t)}_{\L1 (\Omega;\reali)}
    \leq
    K_v \; \norma{w_i}_{\C0 ([0,T];\L1 (\Omega;\reali))} \,.
  \end{equation}
  The regularity $c_i (t) \in \C1 (\Omega;\reali^n)$ directly follows
  from~\ref{it:v}. Concerning the boundedness of $D_x c_i$, observe
  that by~\ref{it:v},
  $\norma{D_x c_i} _{\L\infty ([0,T]\times\Omega;\reali^{n\times n})}
  \leq K_v \; \norma{w_i}_{\C0 ([0,T]; \L1 (\Omega;\reali))}$.

  \paragraph{Step~2: $A_i$ satisfies~\ref{it:H2} in Theorem~\ref{def:hypGen}.}
  Compute:
  \begin{flalign*}
    \norma{A_i}_{\L1([0,T];\L\infty (\Omega;\reali))} = %
    & \int_0^T \esssup_{x \in \Omega} \modulo{\alpha\left(t,x,w_i
        (t)\right)} \d{t}%
    & [\mbox{By~\eqref{eq:18}}]
    \\
    \leq %
    & \int_0^T k_\alpha (t) \left(1 + \norma{w_i (t)}_{\L1
        (\Omega;\reali)}\right) \d{t} %
    & [\mbox{By~\ref{it:alpha}}]
    \\
    \leq %
    & \norma{k_\alpha}_{\L1 (0,T];\reali)}
    \left(1+\norma{w_i}_{\L\infty ([0,T];\L1 (\Omega;\reali))}\right)
  \end{flalign*}
  proving that
  $A_i \in \L1 \left([0,T];\L\infty (\Omega;\reali)\right)$. Moreover,
  \begin{flalign*}
    \int_0^T \tv\left(A_i (t)\right) \d{t} \leq %
    & \int_0^T K_\alpha \left(1 + \norma{\nabla w_i (t)}_{\L1
        (\Omega;\reali^n)}\right) \d{t} %
    & [\mbox{By~\ref{it:alpha}}]
    \\
    = %
    & K_\alpha \; \left(T+\norma{\nabla
        w_i}_{\L1([0,T]\times\Omega;\reali^n)}\right) \,, %
    &[\mbox{By~\ref{item:23} in Proposition~\ref{prop:L2toL1}}]
  \end{flalign*}
  proving that~\ref{it:H2} holds.

  \paragraph{Step~3:
    $B_i \in \L1 \left([0,T]; \L\infty(\Omega;\reali)\right)$.}
  Using~\ref{it:betaFinal} compute:
  \begin{equation}
    \label{eq:26}
    \norma{B_i}_{\L1 ([0,T];\L\infty (\Omega;\reali))} =
    \int_0^T \esssup_{x \in \Omega}
    \modulo{\beta \left(\tau,x,u_i (\tau),w_i (\tau)\right)} \d{\tau}
    \leq \norma{k_\beta}_{\L1 ([0,T];\reali)}
  \end{equation}
  proving Step~3.

  \bigskip

  Hence Theorem~\ref{thm:hypGen} and Theorem~\ref{thm:paraGen} apply,
  ensuring that the sequence $(u_i,w_i)$ can be recursively
  defined. The next steps aim at proving that $(u_i,w_i)$ is a Cauchy
  sequence.

  \paragraph{Step~4: $w_i$ is bounded in
    $\L\infty \left([0,T];\L1 (\Omega;\reali)\right)$ uniformly in
    $i$.}

  By~\eqref{eq:APBound} and~\eqref{eq:26} we have
  \begin{eqnarray*}
    \norma{w_i (t)}_{\L1 (\Omega;\reali)}
    & \leq
    & K \left(\norma{w_o}_{\L1 (\Omega;\reali)} +
      \norma{b}_{\L1 ([0, t] \times \Omega;\reali)} \right)
      \exp\left(K \, \norma{B_{i-1}}_{\L1 ([0, t]; \L\infty (\Omega;\reali))}\right)
    \\
    & \leq
    & K_w
  \end{eqnarray*}
  where we set
  \begin{equation}
    \label{eq:27}
    K_w
    =
    K \left(\norma{w_o}_{\L1 (\Omega;\reali)} +
      \norma{b}_{\L1 ([0, T] \times \Omega;\reali)} \right)
    \exp\left(K \, \norma{k_\beta}_{\L1 ([0, T]; \reali)}\right)
    \,.
  \end{equation}

  \paragraph{Step~5: Bounds uniform in $i$.}

  Note that $c_i (t) \in \C2 (\Omega;\reali^n)$ by~\ref{it:v}
  and~\eqref{eq:18}. Moreover,
  \begin{flalign*}
    \norma{c_i (t)}_{\L\infty (\Omega;\reali^n)} %
    & \leq K_v \norma{w_i}_{\C0 ([0,t];\L1 (\Omega;\reali))} %
    & [\mbox{By~\eqref{eq:28}}]
    \\
    & \leq K_v \, K_w %
    & [\mbox{By~\eqref{eq:27}}]
    \\
    \norma{D_x c_i} _{\L\infty ([0,T]\times\Omega;\reali^{n\times n})}
    & \leq K_v \; \norma{w_i}_{\C0 ([0,T]; \L1 (\Omega;\reali))} %
    & [\mbox{By~\ref{it:v}}]
    \\
    & \leq K_v \, K_w %
    & [\mbox{By~\eqref{eq:27}}]
    \\
    \norma{\nabla \nabla \cdot c_i}_{\L1 ([0,T]\times\Omega;\reali^n)}
    & = \int_0^T \norma{\nabla \nabla \cdot v\left(t,w_i
        (t)\right)}_{\L1 (\Omega;\reali^n)} \d{t} %
    & [\mbox{By~\eqref{eq:18}}]
    \\
    & \leq \int_0^T \norma{D^2_x v\left(t,w_i (t)\right)}_{\L1
      (\Omega;\reali^n)} \d{t}
    \\
    & \leq \int_0^T k_v \left(t,\norma{w_i (t)}_{\L1
        (\Omega;\reali)}\right) \, \norma{w_i (t)}_{\L1
      (\Omega;\reali)} \d{t} %
    & [\mbox{By~\ref{it:v}}]
  \end{flalign*}
  which is bounded since
  $k_v \in \Lloc\infty ([0,T]\times\reali_+;\reali_+)$.

  Using~\eqref{eq:27} and~\ref{it:v}, prepare for later use:
  \begin{eqnarray*}
    \norma{c_{i+1} - c_i}_{\L1 ([0,t];\L\infty (\Omega;\reali^n))}
    & \leq
    & \int_0^t
      \norma{v\left(\tau,w_{i+1} (\tau)\right) - v\left(\tau,w_i (\tau)\right)}_{\L1 (\Omega;\reali^n)} \d\tau
    \\
    & \leq
    & K_v \;
      \norma{w_{i+1}- w_i}_{\L1 ([0,t]\times\Omega;\reali)} \,.
    \\
    \norma{\nabla \cdot(c_{i+1} - c_i)}_{\L1 ([0,t];\L\infty (\Omega;\reali^n))}
    & \leq
    & \int_0^t
      \norma{\nabla \cdot \left(v \left(\tau, w_{i+1} (\tau)\right)
      - v \left(\tau, w_i (\tau)\right)\right)}_{\L\infty (\Omega;\reali)} \d{\tau}
    \\
    & \leq
    & \int_0^t
      k_v\left(\tau, K_w\right) \;
      \norma{w_{i+1} (\tau) - w_i (\tau)}_{\L1 (\Omega;\reali)} \, \d\tau
    \\
    & =
    & \norma{k_v}_{\L\infty ([0,t]\times[0,K_w];\reali)}
      \, \norma{w_{i+1} - w_i}_{\L1 ([0,t]\times\Omega;\reali)}\,.
  \end{eqnarray*}
  Passing to the $A_i$:
  \begin{flalign*}
    & \norma{A_i}_{\L1 ([0,t];\L\infty\Omega;\reali))} %
    \\
    & = \int_0^t \esssup_{x\in\Omega} \modulo{\alpha\left(\tau,x,w_i
        (\tau)\right)} \d{\tau} %
    & [\mbox{By~\eqref{eq:18}}]
    \\
    & \leq \int_0^t k_\alpha (\tau) \left( 1+\norma{w_i (\tau)}_{\L1
        (\Omega;\reali)} \right) \d\tau %
    & [\mbox{By~\ref{it:alpha}}]
    \\
    & \leq \left(1+K_w\right) \norma{k_a}_{\L1 ([0,t];\reali)} \,; %
    & [\mbox{By~\eqref{eq:27}}]
    \\
    & \int_0^t \tv\left(A_i (\tau)\right) \d\tau %
    \\
    & \leq K_a\left(1+\norma{\nabla w_i}_{\L1
        ([0,t]\times\Omega;\reali^n)}\right) %
    & [\mbox{By Step~2}]
    \\
    & \leq K_a+ \frac{K_\alpha}{\mu} \norma{B_{i-1}}_{\L1
      ([0,T];\L\infty (\Omega;\reali))} \; \norma{w_i}_{\L\infty
      ([0,T];\L1 (\Omega;\reali))}
    \\
    & \qquad + \frac{K_\alpha}{\mu} \,
    \norma{b}_{\L1([0,T]\times\Omega;\reali)} + \frac{K_\alpha}{\mu}
    \, \norma{w_o}_{\L1(\Omega;\reali)} %
    & [\mbox{By~\eqref{eq:29}}]
    \\
    & \leq K_a+ \frac{K_\alpha}{\mu} \, K_w \, \norma{k_\beta}_{\L1
      ([0,T];\reali)} %
    & [\mbox{By~\eqref{eq:26}, \eqref{eq:27}}]
    \\
    & \qquad + \frac{K_\alpha}{\mu} \,
    \norma{b}_{\L1([0,T]\times\Omega;\reali)} + \frac{K_\alpha}{\mu}
    \, \norma{w_o}_{\L1(\Omega;\reali)} \,;
    \\
    & \norma{A_{i+1} - A_i}_{\L1 ([0,t];\L\infty (\Omega;\reali))} %
    \\
    & = \int_0^t \norma{A_{i+1} (\tau) - A_i (\tau)}_{\L\infty
      (\Omega;\reali)} \d{\tau}
    \\
    & = \int_0^t \norma{ \alpha\left(\tau,\cdot,w_{i+1} (\tau)\right)
      - \alpha\left(\tau,\cdot,w_i (\tau)\right)}_{\L\infty
      (\Omega;\reali)} \d{\tau}%
    & [\mbox{By~\eqref{eq:18}}]
    \\
    & \leq K_a \, \norma{w_{i+1}-w_i}_{\L1 ([0,t]\times\Omega;\reali)}
    \,. %
    & [\mbox{By~\ref{it:alpha}}]
  \end{flalign*}

  \paragraph{Step~6: $(u_i,w_i)$ is a Cauchy sequence in
    $\C0\left([0,\Delta T]; \L1(\Omega;\reali^2)\right)$ for $\Delta T$ small.}
  Thanks to the above bounds uniform in $i$, all the constants
  appearing in~\ref{item:14} and~\ref{item:15} are bounded uniformly
  by quantities depending on $\Omega$, $u_o$, $a$, $b$ and by the
  constants in the assumptions~\ref{it:omega}, \ref{it:v},
  \ref{it:alpha}, \ref{it:a}, \ref{it:betaFinal}, \ref{it:b}. Hence,
  \ref{item:14} and~\ref{item:15} yield
  \begin{eqnarray*}
    &
    & \norma{u_{i+1} (t) - u_i (t)}_{\L1 (\Omega;\reali)}
    \\
    & \leq
    & \O
      \left(
      \norma{c_i - c_{i-1}}_{\L1 ([0,t];\L\infty (\Omega;\reali^n))}
      +
      \norma{\nabla \cdot(c_i - c_{i-1}}_{\L1 ([0,t];\L\infty (\Omega;\reali^n))}
      \right)
    \\
    &
    & \qquad +
      \O
      \left(
      \norma{A_i-A_{i-1}}_{\L1 ([0,t];\L\infty (\Omega;\reali))}
      \right)
    \\
    & \leq
    & \O \; \norma{w_i - w_{i-1}}_{\L1 ([0,t]\times\Omega;\reali)} \,.
  \end{eqnarray*}
  To compute the distance
  $\norma{w_{i+1} (t) - w_i (t)}_{\L1 (\Omega;\reali)}$,
  apply~\ref{item:10} thanks to~Step~3 and get
  \begin{flalign*}
    &\norma{w_{i+1} (t) - w_i (t)}_{\L1 (\Omega;\reali)}
    \\
    \leq %
    & K^2 \left( \norma{w_o}_{\L1 (\Omega;\reali)} + \norma{b}_{\L1
        ([0,t]\times\Omega;\reali)} \right) e^{K (\norma{B_i}_{\L1
        ([0,t];\L\infty (\Omega;\reali))} + \norma{B_{i-1}}_{\L1
        ([0,t];\L\infty (\Omega;\reali))})}
    \\
    & \quad \times \norma{B_i - B_{i-1}}_{\L1 ([0,t];\L\infty
      (\Omega;\reali))} %
    & [\mbox{By~\ref{item:10}}]
    \\
    \leq %
    & K^2 \, K_\beta \left( \norma{w_o}_{\L1 (\Omega;\reali)} +
      \norma{b}_{\L1 ([0,t]\times\Omega;\reali)} \right) e^{2K
      \norma{k_\beta}_{\L1 ([0,t];\reali)}} %
    & [\mbox{By~\eqref{eq:26}}]
    \\
    & \quad \times \left( \norma{u_i - u_{i-1}}_{\L1
        ([0,t]\times\Omega;\reali)} %
      + %
      \norma{w_i - w_{i-1}}_{\L1([0,t]\times\Omega;\reali)} \right)\,. %
    & [\mbox{By~\ref{it:betaFinal}}]
  \end{flalign*}
  As a consequence,
  \begin{eqnarray*}
    &
    & \norma{u_{i+1} - u_i }_{\C0([0,\Delta T]; \L1(\Omega;\reali))}
      +
      \norma{w_{i+1} - w_i }_{\C0([0,\Delta T]; \L1(\Omega;\reali))}
    \\
    & \leq
    & \O \left(
      \norma{u_i - u_{i-1} }_{\L1([0,\Delta T]\times\Omega;\reali)}
      +
      \norma{w_i - w_{i-1} }_{\L1 ([0,\Delta T]\times\Omega;\reali)}
      \right)
    \\
    & \leq
    & \O \, \Delta T \left(
      \norma{u_i - u_{i-1} }_{\C0([0,\Delta T]; \L1(\Omega;\reali))}
      +
      \norma{w_i - w_{i-1} }_{\C0([0,\Delta T]; \L1(\Omega;\reali))}
      \right)
  \end{eqnarray*}
  proving that for $\Delta T$ small, $(u_i,w_i)$ is a Cauchy sequence
  in $\C0 \left([0,\Delta T];\L1 (\Omega;\reali)\right)$. Call $(u,w)$
  the corresponding limit.

  \paragraph{Step~7: Problem~\eqref{eq:1} admits a global solution in
    the sense of Definition~\ref{def:def_sol}.}

  Note that by~\ref{it:v}, \ref{it:alpha} and~\ref{it:betaFinal}, we
  can pass to the $\C0 \left([0,\Delta T];\L1 (\Omega;\reali)\right)$
  limit also in~\eqref{eq:18}.  Hence, the Dominated Convergence
  Theorem allows to pass to the limit in~\eqref{eq:23}, so that
  Definition~\ref{def:hypGen} applies to $u$. A further use of
  Lemma~\ref{lem:weakComplete} allows to pass to the limit also
  in~\eqref{eq:22}, proving the existence of a solution on the time
  interval $[0,\Delta T]$.

  Further iterations of the above procedure yield a solution, say
  $(u,w)$, to~\eqref{eq:1} in the sense of
  Definition~\ref{def:def_sol}. Call $[0,T_*\mathclose[$, for a
  $T_* > \Delta T$, the biggest time interval on which $(u,w)$ can be
  extended.

  Define $B (t,x) = \beta\left(t,x,u (t), w
    (t)\right)$. By~\ref{it:betaFinal},
  $B \in \L1\left([0,T_*]; \L\infty (\Omega;\reali)\right)$. Hence,
  Theorem~\ref{thm:paraGen} ensures that problem~\eqref{eq:13} admits
  a solution $w$ on $[0, T_*]$ in the sense of
  Definition~\ref{def:paraGen}. By~\ref{item:17},
  $w \in \C0 \left([0,T_*];\L1 (\Omega;\reali)\right)$.

  Define $c (t,x) = v \left(t,w (t)\right) (x)$,
  $A (t,x) = \alpha \left(t,x,w (t)\right)$ and repeat the same
  computations as in Step~1 and Step~2 to obtain that
  problem~\eqref{eq:21} admits as solution $u$ on $[0,T_*]$ in the
  sense of Definition~\ref{def:hypGen}. By~\ref{item:16}, $u$ is
  continuous on $[0,T_*]$.

  Thus, $(u,w)$ is extended up to time $T_*$. If $T_* <T$, the above
  procedure can be repeated with reference to problem~\eqref{eq:1}
  with initial datum $(u,w) (T_*)$ assigned at time $T_*$, obtaining
  an extension of $(u,w) $ beyond time $T_*$, hence contradicting the
  maximality of $T_*$, unless $T_* = T$. This proves the global
  existence stated in~\ref{item:27}.

  The continuity at~\ref{item:31} follows from~\ref{item:16}
  and~\ref{item:17}.

  The continuous dependence and stability estimates~\ref{item:28},
  and~\ref{item:29} follow through long and tedious computations based
  on the estimates obtained so far. More precisely, to
  prove~\ref{item:28}, use~\ref{item:13}, \ref{item:14},
  \ref{item:15}, \ref{item:9}, \ref{item:10} and repeat the
  computations in Step~3, 4, 5 and~6. To prove~\ref{item:29} the
  procedure is entirely similar, also using~\ref{item:32}.

  The proof of~\ref{item:33} is a direct consequence of~\ref{item:12}
  and~\ref{item:8}.
\end{proofof}

\paragraph{Acknowledgment:} The authors acknowledge the PRIN 2022
project \emph{Modeling, Control and Games through Partial Differential
  Equations} (D53D23005620006), funded by the European Union - Next
Generation EU. The first two authors also acknowledge the 2024~GNAMPA
project \emph{Modelling and Analysis through Conservation Laws}.

\paragraph{Data Availability Statement:} Data sharing not applicable
to this article as no data set were generated nor analysed during the
current study.

{ \small

  \bibliography{neumann}

  \bibliographystyle{abbrv}

}

\end{document}